\documentclass[12pt,twoside]{article}

\usepackage{latexsym, amssymb, amsmath, theorem}
\usepackage{epsfig}

\numberwithin{equation}{section}

\theoremstyle{change}
\theorembodyfont{\itshape}
\newtheorem{theorem}{Theorem}[section]
\newtheorem{proposition}[theorem]{Proposition}
\newtheorem{lemma}[theorem]{Lemma}
\newtheorem{corollary}[theorem]{Corollary}

\theorembodyfont{\rmfamily}
\newtheorem{definition}[theorem]{Definition}

\newtheorem{example}[theorem]{Example}
\newtheorem{remark}[theorem]{Remark}
\newenvironment{proof}{{\noindent \textbf{Proof}\,\,}}{\hspace*{\fill}$\Box$\medskip}

\newcounter{amarker}
\renewcommand{\theamarker}{\textup{(\arabic{amarker})}}

\newcounter{rmarker}

\makeatletter
\renewcommand{\@makecaption}[1]{%
\begin{center}#1\end{center}%
}
\makeatother
\def\ccd{\mathbb C^2}
\title{Upper bounds of topology of complex polynomials in two
variables\protect\footnote{AMS Classification numbers 14D05,
30C15, 30C10}}
\author{A.A.Glutsyuk\protect\footnote{Present Institution:
Poncelet Laboratory
(UMI 2615 of CNRS and Independent University of
Moscow)}\protect\footnote{Permanent address:
CNRS, UMR 5669, Unit\'e de Math\'ematiques pures
et  appliqu\'ees, M.R.,
\'Ecole Normale Sup\'erieure de Lyon, 46, all\'ee d'Italie, 69364 Lyon Cedex 07
France}}
\begin{document}
\maketitle

\hskip5cm {\it Dedicated to Mikhail Anatolievich Tsfasman}

\hskip5cm {\it on the occassion of his 50-th birthday}

\renewcommand\d{\partial}
\renewcommand\c{\mathbb C}
\def\buf{B_{K,U} (f)}
\newcommand\bui{B_{K,U} (I)}
\newcommand\e{\varepsilon}
\newcommand\ph{\varphi}
\newcommand\nbd{neighborhood}
\tableofcontents

\def\wt#1{\widetilde#1}
\def\tb{\tilde{\beta}}
\def\dx{\operatorname{dx}}
\def\dy{\operatorname{dy}}
\def\diag{\operatorname{diag}}
\def\grad{\operatorname{grad}}
\def\De{\delta}
\def\oi{\omega_i}
\def\oj{\omega_j}
\def\var{\varepsilon}
\def\gt{\gamma_t}
\def\gto{\gamma_{t_0}}
\renewcommand\o{\omega}
\def\cc{\mathbb C}
\def\an{a_1,\dots,a_{n^2}}
\def\can{\cc\setminus\{\an\}}
\def\guk{\gamma(U,K)}
\def\a{\alpha}
\def\gi{\delta_i}
\def\gj{\delta_j}
\newcommand\D{\Delta}
\def\zz{\mathbb Z}
\def\mbi{\mathbb I}
\def\dj{\delta_j}

\section{Introduction, main results and the plan of the paper}

\begin{definition} \label{defgen} We say that a homogeneous
polynomial is {\it generic}, if it has only simple
zero lines.
\end{definition}

Let $H(x,y)$ be a complex polynomial in two variables,
$h$ be its highest homogeneous part, $H'=H-h$ be its
lower terms. We assume that

- $h$ is generic,

- $0$ is a critical point of $H$,

$$H'-H(0)\not\equiv0.$$

It is well-known that under these assumptions
 $H$ has at least two distinct critical values.
In the present paper we prove quantitative versions of this
statement (Theorems \ref{topxy} and \ref{t3.5}).
Assuming that $H$ is appropriately normalized by affine variable
changes in the image and the preimage
(see Definition \ref{defnor}), we give explicit
upper bounds of the following quantities:

- the lower terms $H'=H-h$ (Addendum to Theorem \ref{topxy});

- the minimal size of a closed bidisc containing a
deformation retract
(i.e., all the nontrivial topology) of a level curve
$$S_t=\{ H=t\} \ \text{(Theorem \ref{topxy})};$$

- the minimal lengths of representatives of cycles in $H_1(S_t,\zz)$
vanishing along appropriate paths from $t$ to the critical values
(Lemma \ref{l(n)} and Corollary \ref{cl(n)} in 1.5,
these representatives are contained in the previous bidisc);

- the intersection indices of the latter cycles
(Theorem \ref{tindex} in 1.6).

The above-mentioned upper bounds
of locations and lengths of vanishing cycles were used in \cite{gi}
(a joint paper with Yu.S.Ilyashenko),
where we studied zeros of Abelian integrals of polynomial
1- forms over real ovals of real polynomials. In \cite{gi} we gave
an upper bound of the number of zeros for a wide class of
real polynomials of arbitrary degree $n+1\geq3$ and any polynomial
1- form of smaller degree. This class consists of ultra-Morse
polynomials $H$ (see Definition \ref{dump}). The estimate
of the number of zeros
is done in terms of the minimal gap between
two zero lines of the higher homogeneous part $h$
and the ratio of that between  two critical values of $H$
over the diameter of the critical value set.
The key object considered in the proofs in \cite{gi} is the
Abelian integral matrix
\begin{equation} \mbi(t)=(I_{ij}(t)), \
I_{ij}(t)=\int_{\dj\subset S_t}\oi,\ \text{where}\label{2.3}
\end{equation}
$\oi$, $i=1,\dots,\mu=n^2$, are fixed monomial forms
of appropriate degrees, $\dj\subset H_1(S_t,\zz)$
are appropriate cycles (marked vanishing cycles) that
form a basis.

\begin{remark} The matrix elements $I_{ij}(t)$ are multivalued
holomorphic functions with branchings at the critical values of
$H$.
\end{remark}

In the present paper we prove the following statements used
in \cite{gi}:

- upper bound of the elements of the Abelian integral matrix
(\ref{2.3}) (see 1.7);

- lower bound of its determinant
\begin{equation}\Delta(t)=\D_{H,\Omega}(t)=\det\mbi(t),\
\Omega=(\omega_1,\dots,\omega_{n^2})\label{dett}
\end{equation}
(see 1.8), which is called the {\it period determinant}.
We prove these statements in full generality, without assumption
that $H$ is real.

The above-mentioned upper bound of $|I_{ij}(t)|$ is implied
by the previously-mentioned upper bound of
locations and  lengths of vanishing cycles.
The upper bound of the latter lengths is proved by using
the upper bound of topology of level curve and Bezout's theorem.

It is well-known that if $\omega_i$
are  homogeneous polynomial 1- forms
of appropriate degrees, then the period determinant is a polynomial
of the type

\begin{equation}\D(t)=C(h,\Omega)\prod_{i=1}^{n^2}(t-a_i),
\ \text{where}
\label{2.20}\end{equation}
$a_i$ are the critical values of $H$ (\cite{ga}, see \cite{g} for more
details, where an explicit formula for the constant
$C(h,\Omega)$ was obtained). The proof of the lower bound of
$\D(t)$   is based on the latter formula.

The upper bound of topology and the lower bound of the period
determinant are the principal results of the paper: their
proofs take its most part.

\subsection{The plan of the paper} The results concerning
upper bounds of topology of level curve and of lower terms of $H$
are stated in 1.2 and proved in 1.2 and Section 2.

The definition of vanishing cycles is recalled in 1.3.

Upper bounds of lengths of their
appropriate (canonical) representatives
are stated and proved in 1.5.
These canonical representatives are defined in 1.4;
their projections to the $x$- axis are piecewise-linear curves.
The number of pieces is estimated at the same place (Proposition
\ref{pnopieces} and Corollary \ref{nopieces}). The
intersection indices of vanishing cycles are estimated in 1.6
(using Corollary \ref{nopieces}). Proposition \ref{pnopieces}
is the key statement used
in the estimates of lengths and intersection indices.

Upper bounds (used in \cite{gi})
of the integrals $I_{ij}(t)$ are stated and proved in 1.7.

The lower bound of the period determinant is stated in 1.8 and
proved in Section 3.

\subsection{Upper bounds of topology and lower terms}
To state the results from the title of the Subsection,
we have to normalize the polynomial in appropriate way by
affine variable changes. To specify the normalization,
let us firstly introduce the following
definitions (see also \cite{gi}, where some of them were
also introduced).

\begin{definition} \label{d1.6} The
{\it norm} of a homogeneous polynomial is
the maximal value of its module on the unit sphere; this norm
is denoted by $\| h \|_{\max}.$ The norm of a nonhomogeneous
polynomial is the sum of norms of its homogeneous parts.
\end{definition}

For any $r>0$ denote
$$D_r=\{z\in\cc\ | \ |z|<r\}, \ D_r(a)=\{ z\in\cc\ | \
|z-a|<r\}.$$
For any $X,Y>0$ denote
$$D_{X,Y}=\{(z,w)\in\ccd\ | \ |z|<X,\ |w|<Y\}.$$

\begin{definition} \label{defnor}
A polynomial $H$ as at the beginning of the
paper is said to be {\it weakly normalized}, if $||h||_{max}=1$ and
all its critical values are contained in the closed disc
$\overline D_2$. It is said to be {\it normalized}, if
the previous conditions hold and
there is no smaller disc containing the critical values.
\end{definition}

\begin{remark} Any polynomial $H$ as at the beginning of the
paper can be transformed to a normalized one by
(nonunique) affine variable
changes in the image and the preimage so that the highest
homogeneous part remains unchanged up to multiplication by
constant. (The previous definition is a slightly extended version
of an analogous one  from \cite{gi}.)
\end{remark}

Let us introduce the following function on the space of
normalized polynomials.

\begin{definition} \label{d1.3} For any polynomial
$H$ of degree $n+1$ with a generic highest homogeneous part $h$
put $c_1(H)$ to be
$n$ multiplied by the smallest distance between two lines in the
zero locus of $h$. The distance between
two lines is taken in sense of Fubini-Study metric on the
projective line $\mathbb CP^1.$ Let
$$c'(H) = \min (c_1(H),1).$$
\end{definition}

\begin{definition} \label{deftop}
We say that the topology of a level curve
$S_t=\{ H=t\}$ is contained in a bidisc $D_{X,Y}$, if
 the difference $S_t \setminus D_{X,Y}$ consists of
$n +1= \deg H$ punctured topological disks, and the restriction of
the projection $(x,y) \mapsto x$ to any of these disks is a
biholomorphic map onto $\{ x \in \mathbb C\ |\ X < |x| < \infty \} .$
\end{definition}

\begin{theorem} \label{topxy} For any weakly
normalized polynomial $H$ of degree $n+1\geq3$
the Hermitian basis in $\mathbb C^2$ may be so chosen that
the  topology
of all the level curves $S_t$ for $|t| \le 5$ will be located in a
bidisc $D_{X,Y}$ with
\begin{equation}
X\leq Y \le  {(c'(H))}^{-14n^3}n^{65n^3}=R_0.
 \label{xy<<}\end{equation}
The choice of basis in $\ccd$ is specified below.
\end{theorem}

{\bf Addendum} In the conditions of the Theorem the norm of
the lower terms $H'=H-h$ admits the following upper bound:
\begin{equation}||H'||_{max}<(c')^{-13n^4}n^{64n^4}.
\label{nonlin}\end{equation}

\medskip

The Theorem and the Addendum are proved in Section 2.

\begin{corollary} \label{topxyt}
In the conditions of the previous Theorem for any
$t\in\cc$ the topology of the level curve $S_t$ is contained in
the bidisc $D_{X,Y}$ with
\begin{equation}
X\leq Y\leq R_0\chi(|t|),
\chi(t)=\max\{1,(\frac{|t|}5)^{\frac1{n+1}}\}.\label{defchi}
\end{equation}
\end{corollary}

\begin{proof} If $|t|\leq 5$, then the statement of the Corollary
follows immediately from the Theorem. If $|t|>5$, consider
the rescaled
polynomial $\wt H(x,y)=\frac5{|t|}H(|\frac t5|^{\frac1{n+1}}x,
|\frac t5|^{\frac1{n+1}}y)$. It is also weakly normalized.
Indeed, its highest homogeneous part coincides with that of $H$.
Its critical values are equal to
$\frac5{|t|}<1$ times those of $H$, thus, their modules
are no greater than 2. This implies weak normalizedness.
The level curve $S_t=\{ H=t\}$ is transformed by the previous
rescaling to the level curve
$\{\wt H=t\frac5{|t|}\}$, which corresponds to a level value
of $\wt H$ with module 5 by construction.
Hence, by the Theorem applied to $\wt H$,
the topology of the latter curve
is contained in the bidisc
$D_{X',Y'}$, $X'\leq Y'\leq R_0$. Applying the inverse rescaling
yields that the topology of the curve $S_t$ is contained in
$D_{X,Y}$, $X\leq Y\leq R_0|\frac t5|^{\frac1{n+1}}$. This proves
the Corollary.
\end{proof}

{\bf Choice of basis in $\ccd$.} Theorem \ref{topxy} will be proved
for orthogonal coordinates $(x,y)$ in $\mathbb C^2$
satisfying the
following inequality (their existence is implied by the next
Proposition):
\begin{equation}\text{the distance
 of the}\ y- \ \text{axis to the zero lines of}\ h\ \text{is greater  than}\
 \frac1{\sqrt n}.\label{disty0}\end{equation}

\begin{proposition} \label{pdisty0} For any $m\in\mathbb N$ and
any tuple of $m+1$ complex lines in $\mathbb C^2$ passing
through the origin there exists another complex line through the
origin whose distance to each line of the tuple is
greater than $\frac1{\sqrt m}.$
\end{proposition}

\begin{proof} Take the $\frac1{\sqrt m}$ neighborhood in the sphere
$\mathbb{CP}^1$ of each line of the tuple.
The area of the neighborhood is less than that of an Euclidean disc
of the same radius, i.e., less than $\pi\frac1m$. The
union of the latter neighborhoods through all the lines from
the tuple
does not cover the whole sphere  (and this implies the Proposition).
Indeed, the area of the union is less than
$\frac{m+1}m\pi<2\pi$, which is less
than the area $4\pi$ of the whole sphere. The Proposition is proved.
\end{proof}

\begin{definition} \label{dunitsc} A complex polynomial $H=h+H'$
with a generic highest homogeneous part $h$ is said to be {\it unit-scaled}, if
$$||h||_{max}=1\geq ||H'||_{max}.$$
\end{definition}

Let us give a brief description of proof of Theorem \ref{topxy}.
This Theorem and its Addendum hold true automatically
for unit-scaled polynomials $H$, and a much stronger bound of topology follows
almost immediately by elementary inequalities. Namely,
one  uses the fact that outside a ball whose radius can be
explicitly  estimated from above the foliation $H=const$
resembles the foliation by level curves of the homogeneous
polynomial $h$. This is done in Subsection 2.3. The results
of 2.3 imply the following

\begin{proposition} \label{phl} For any unit-scaled polynomial
$H$ of degree $n+1\geq3$
and any $t\in\cc$, $|t|\leq5$, the topology of the level curve
$S_t$ is contained in the bidisc $D_{X_n,Y_n}$,
\begin{equation}X_n=n^{7n}(c'(H))^{-2n},\ \
Y_n=n^nX_n=n^{8n}(c'(H))^{-2n}.
\label{xnyn}
\end{equation}
\end{proposition}

In the general case, when the norm of lower terms may be large,
for the proof of Theorem \ref{topxy}, we change the normalization
of $H$ to make the norm of lower terms unit, and we prove an upper
bound of the rescaling rate. This bound is based on the
next Theorem, which is a quantitative analogue of the first
statement mentioned at the beginning of the paper. Roughly
speaking, it says that
if the difference $H'-H(0)$ is not too small, then
the maximal distance between critical values of $H$
admits an explicit lower bound. To formulate  it, we change
normalization again in such a way that in addition $H(0)=0$.

\def\crp{centrally-rescaled }

\begin{definition} A polynomial $H$ as at the beginning of the paper
is said to be {\it \crp}, if $H(0)=0$ and $||h||_{max}=||H'||_{max}=1$.
\end{definition}

\begin{remark} Each polynomial $H$ from the beginning of the paper
can be transformed to a \crp one
by homothety in the preimage and an affine
transformation in the image.
\end{remark}

\begin{theorem} \label{t3.5} Each \crp polynomial $H$ of degree
$n+1\geq3$ has at least one
critical value with module no less than
\begin{equation}\delta_0=(c'(H))^{13n^4}n^{-63n^4}.\label{fdel}\end{equation}
\end{theorem}

The proof of Theorem \ref{t3.5} takes the most part of Section 2.
Theorem \ref{topxy} will
be deduced from Theorem \ref{t3.5} and Proposition \ref{phl}
in 2.2.

For the proof of Theorem \ref{t3.5} we consider projections
$\pi_t:S_t\to Ox$ of level curves to the $x$- axis. The proof is done by
combining topological arguments and analysis of arrangement
configuration of critical values of projections.
The topological arguments are based on connectivity
of the intersection graph of marked basic vanishing cycles and
Picard-Lefschetz theorem [1]. The main technical result
of the analysis of critical values of projection is a lower bound of
the maximal module of a critical value of $\pi_0$
(Lemma \ref{distcr0} stated in 2.1).
The topological and analytical
arguments are related to each other via studying
appropriate (canonical) representatives of the vanishing cycles
whose projections are piecewise-linear curves with vertices
at the critical values of projection.

\subsection{Vanishing cycles}
All the definitions and the statements of the present
Subsection are contained in \cite{1} and \cite{gi}.

The basic vanishing cycles are usually defined in homology
of level curves of an ultra-Morse polynomial, see the following
Definition.

\begin{definition} \label{dump}
A complex polynomial of degree $n+1$
with a generic highest homogeneous part
is said to be {\it ultra-Morse}, if all its
critical values are simple (then their number is equal to $\mu=n^2$).
\end{definition}

Firstly we recall the  definition of a local vanishing cycle.

\begin{lemma} {\bf (Morse lemma).}
A holomorphic function having a Morse
critical point may be transformed to a sum of a nondegenerate
quadratic form and a constant term by an analytic change of
coordinates near this point.
\end{lemma}

\begin{corollary} Consider a holomorphic function in $\mathbb C^2$
having a Morse critical point with a critical value $a.$ There exists a ball
centered at the critical point whose intersection with
each level curve corresponding to a value close to $a$ of the
function is diffeomorphic to an annulus.
\end{corollary}

\begin{definition} A generator of the first homology
group of the latter intersection annulus (considered as a cycle in
the homology of the global level curve) is called a
{\it local vanishing cycle} corresponding to $a.$
\end{definition}
A local vanishing cycle is well defined up to change of
orientation.
\begin{definition} \label{dreg}
Let $H$ be a ultra-Morse polynomial, $a_j$ be its critical values,
$j=1,\dots,\mu$.
A path $\alpha_j: [0,1] \to \mathbb C$ is called {\it regular}
provided that
\begin{equation}
\alpha_j(1)=a_j, \ \alpha_j [0,1)
\ \text{contains no critical value of}\ H.
\label{2.1}
\end{equation}
\end{definition}
\begin{definition} \label{d2.3} Let $\alpha_j$ be a regular path,
$t_0=\a_j(0)$, $s \in
[0,1]$ be close to $1, \ \delta_j(t), \ t = \alpha_j(s),$ be a
local vanishing cycle on $S_t$ corresponding to $\a_j.$ Consider
the extension of $\delta_j$ along the path $\a_j$ up to a
continuous family of cycles $\delta_j(s)$ in complex level curves
$H = \alpha_j(s).$ The homology class
$\delta_j = \delta_j(0)\in H_1(S_{t_0},\mathbb Z)$ is
called a {\it cycle vanishing along $\a_j.$}
\end{definition}
\begin{definition} \label{d2.4} Consider a set of regular paths $\alpha_1, \dots , \alpha_\mu,$
see (\ref{2.1}), with a common starting point $t_0$.
Suppose that these paths are not pairwise and self
intersected. Then the set of cycles $\delta_j \in H_1(S_{t_0},
\mathbb Z)$ vanishing along
$\alpha_j, \ j = 1. \dots , \mu ,$ is called a
{\it marked set of vanishing cycles} on the level curve $H = t_0.$
\end{definition}

\begin{lemma}\label{lbas} Any marked set of vanishing
cycles is a basis in the first integer homology group of the level curve.
\end{lemma}

The extension of a cycle in $H_1(S_t,\zz)$
along a loop (avoiding critical values)
defines an operator from the homology
group to itself called {\it monodromy operator}.

\begin{lemma}\label{lmon} The images of any vanishing cycle
under the monodromy operators
along all the loops generate the homology group.
\end{lemma}

The Lemma follows from Picard-Lefschetz theorem and the
connectedness of the intersection graph of marked set of vanishing
cycles \cite{1}.

\subsection{Canonical representatives of vanishing cycles}

Below we define the representatives of vanishing cycles for which we
prove upper bounds of lengths and locations. To do this, we use
the following properties of projection of local level curve
near critical point, which are well-known and follow immediately
from definition.

\def\var{\varepsilon}

\begin{proposition} \label{doublecr}
Let $H(x,y)$ be an analytic function having
a Morse critical point at 0, $H(0)=0$. The singular level curve
$H=0$ is locally a union of two regularly embedded analytic
curves intersecting at 0 transversally (denote $L_i$, $i=1,2$,
their tangent lines at 0). Let $\pi_t:S_t\to\mathbb C$ be a
projection along a line transversal to both $L_i$. There exist
a neighborhood $U$ of the critical point and an $\var>0$
such that for any $t\neq0$ close to 0 the restriction
$\pi_t|_{S_t\cap U}$ defines a degree two branched covering
$$\pi_t:\pi_t^{-1}(D_{\var})\cap U\to D_{\var}.$$
\end{proposition}

{\bf Addendum.} The latter covering has
two critical values (branching points) $x_1(t)$, $x_2(t)$
confluenting to 0, as $t\to0$. The union of the
two liftings of the segment $[x_1,x_2]$ under the
covering is a closed curve representing a local vanishing cycle
(the liftings are taken with inverse orientations).

\begin{definition} Let in the previous Proposition the projection
$\pi_t$ be made along the $y$- axis. Then
the corresponding closed curve from the Addendum is called
the {\it canonical representative} of the local vanishing cycle
(see Fig.1).
\end{definition}

\begin{figure}[ht]
  \begin{center}
   \epsfig{file=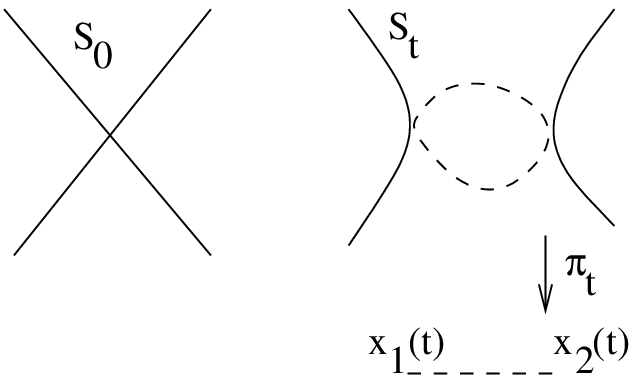}
    \caption{}
    \label{fig:1}
  \end{center}
\end{figure}

Everywhere below (whenever the contrary is not specified) we assume that
$H$ is a ultra-Morse polynomial and the orthogonal coordinates $(x,y)$
are chosen so that

\begin{equation}
\text{the zero lines of}\ h\ \text{and
the tangent lines of the critical level curve branches of}\ H
\label{locbr}\end{equation}

at the critical points are transversal to the $y$- axis.

Then the canonical representatives of the local vanishing
cycles are well-defined.

\begin{remark} One can always chose orthogonal coordinates in $\ccd$
so that both statements (\ref{disty0}) and (\ref{locbr}) hold true.
Indeed, if (\ref{disty0}) holds true (which is an open condition),
then zero lines of $h$ are automatically transversal to the
$y$- axis. Its transversality to the tangent lines of the critical
level curves at the critical points can be achieved by
small perturbation.
\end{remark}

For any noncritical value $t$ denote
$$CS_t=\{\text{critical values of the projection}\ \pi_t:S_t\to Ox\}.
$$

Now we construct the canonical representatives of the global
vanishing cycles along appropriate piecewise-linear paths,
see the next Definition. To do this, we use the following
\begin{remark} \label{reta}
There exists a $\mu'\in\mathbb N$ (depending on $H$)
such that for all
but finite number of values of $t$ the cardinality of the
set $CS_t$ (without multiplicities)
is equal to $\mu'$ (algebraicity). For a typical
polynomial $H$ all the critical values of $\pi_t$
are simple and $\mu'=n(n+1):=\eta(n)$ (Bezout's theorem).
In general, some critical values of the projection
may be multiple for all $t$
(i.e.,  $\mu'<\eta(n)$). This is the case for the
homogeneous polynomial $x^{n+1}+y^{n+1}$, where $\mu'=n+1$.
\end{remark}

\begin{definition} \label{dcreg}
A piecewise-linear regular (see Definition \ref{dreg})
path $\a:[0,1]\to\cc$ with a finite number of edges
is said to be {\it critically-regular}, if for any $t\in\a[0,1)$
the cardinality of the set $CS_t$ is maximal (i.e., equal to
$\mu'$).
\end{definition}

\def\crpa{\text{critically-regular} }

\begin{remark} \label{rkcreg} Any regular path is homotopic
(outside the critical values of $H$) to some \crpa path.
\end{remark}

Let $a\in\cc$ be a critical value of $H$,
$A\in\ccd$ be the corresponding critical point,
$x(A)$ be its $x$- coordinate. Let
$t_0\in\cc$ be a noncritical value of $H$, $\a:[0,1]\to\cc$
be a \crpa path from $t_0$ to  $a$. Let
$$CS_t=\{ x_1(t),\dots,x_{\mu'}(t)\}.$$
The functions $x_i(t)$ are multivalued, but their restriction
pull-backs $x_i(\a(\tau))$ along $\a$ have well-defined continuous
branches with disjoint graphs (critical regularity).
Let $x_1(\a(\tau))$, $x_2(\a(\tau))$ be those of
them that confluent to $x(A)$, as $\tau\to1$, $\wt\delta(\tau)$
be the canonical representative of the corresponding
local vanishing cycle. This representative is well-defined
for $\tau$ close to 1 as  a double lifting to $S_t$
of the segment $[x_1(t),x_2(t)]$, $t=\a(\tau)$. Let us construct
its continuous extension
(as a family of closed curves in $S_{\a(\tau)}$)
for all $\tau\in[0,1)$. Then the curve $\wt\delta(0)$ represents the
global vanishing cycle along $\alpha$ and will be called its
{\it canonical representative}.

Case 1: $x_i(\a(\tau))\notin(x_1(\a(\tau)),x_2(\a(\tau)))$
whenever $\tau\in[0,1)$, $i\neq1,2$. Then the previous
double lifting $\wt\delta(\tau)$ of $[x_1,x_2]$
is well-defined and continuous for all $\tau\in[0,1)$.

Case 2: $\alpha$ is a segment and
there exists a parameter value $\tau_1$
for which
\begin{equation} \text{some}\ x_i,\ i\neq1,2,\ \text{meet}\
(x_1,x_2), \ \text{denote}\ N(\tau_1)\ \text{their number}.
\label{ntl}\end{equation}
(No $x_i$, $i\neq1,2$, can meet $x_1$, $x_2$, by
critical regularity of $\a$.) For $\tau=\tau_1$ the previous
points $x_i$ split
$[x_1,x_2]$ into a union of $N(\tau_1)+1$ smaller segments
$[x_{i_r},x_{i_{r+1}}]$, $r=0,\dots,N(\tau_1)$.
The new segments (whose ends are functions in $\tau$) do not
meet other $x_j$ whenever $\tau$ is close enough to $\tau_1$
(by construction and critical regularity). Then the family
$\wt\delta(\tau)$ extends continuously
to the values $\tau\leq\tau_1$ close to $\tau_1$
as the union of appropriate
double liftings of the previous segments $[x_{i_r},x_{i_{r+1}}]$.
This extension is well-defined until we cross a parameter
value $\tau_2<\tau_1$ for which some of the intervals
$(x_{i_r},x_{i_{r+1}})$
(which are the edges of the projection of $\wt\delta(\tau)$)
meet some other $x_j$' s. Then we repeat the previous extension
of $\wt\delta$ to smaller values $\tau<\tau_2$, etc.
By algebraicity, after a finite number of steps we will extend
$\wt\delta(\tau)$ to all the values $\tau\in[0,1)$
(the number of steps and edges of the projection
is estimated below using Bezout's theorem).

Case 3 (general): the path $\a$ has several edges. We extend
$\wt\delta$ by induction in the number of edges of $\a$. The
induction base (one edge) is given by the previous construction
(Case 2).

Induction step. Let $\a[0,\tau']$ be the first edge of $\a$.
By the induction hypothesis, the family of closed curves
$\wt\delta(\tau)$ is extended to $\tau\in[\tau',1)$.
The previous construction
applied to each edge of the projection of $\wt\delta(\tau')$
(instead of $[x_1,x_2]$) yields the desired extension of
$\wt\delta(\tau)$ to $\tau\in[0,\tau')$.

\begin{definition} Thus constructed closed curve $\wt\delta(0)$
is called the {\it canonical representative} of the
cycle vanishing along $\a$.
\end{definition}

\begin{proposition} \label{pieces} The canonical representative
in a level curve $S_t$ of a cycle
vanishing along a \crpa path
is projected onto a piecewise-linear curve in the
complex $x$- axis
with vertices in $CS_t$ (the edges are minimal:
no edge interval contains a point of $CS_t$).
More precisely, the representative is a finite
union of couples of arcs; the arcs from each couple
are disjoint (maybe except for their ends) and are projected
bijectively onto the same edge.
If the topology of $S_t$
is contained in a bidisc $D_{X,Y}$, then the
canonical representative is contained in the same bidisc.
For any given \crpa path
$\a:[0,1]\to\cc$ the canonical representatives
of cycles vanishing along smaller paths $\a|_{[\tau,1]}$,
$0\leq\tau<1$, depend continuously on $\tau$.
\end{proposition}

This Proposition follows immediately from construction. Its
statement concerning the location of the representative in
$D_{X,Y}$ follows from the fact that all the critical values
of projection should lie in the disc $D_X$ (by the definition of
the latter) and from the convexity of this disc.

\begin{remark} The projection of a canonical representative of
vanishing cycle is a piecewise linear curve with a given order
of pieces (up to inversion). It may happen that several distinct
arc couples are projected onto one and the same edge
(the number of all the arc couples over a given edge is called
its {\it multiplicity}).
The number of arc couples from the previous Proposition
is equal to the number of projection edges with multiplicities.
\end{remark}

The upper bounds of lengths and intersection indices of vanishing
cycles are based on the next Proposition (proved below),
which is proved by using Bezout's theorem.

For any piecewise-linear curve $\Gamma$ denote
$$m(\Gamma)=\#(\text{edges of}\ \Gamma).$$
The number of arc couples of a canonical representative
$\wt{\delta}$
of vanishing cycle (or equivalently, the number of edges
(with multiplicities) of its projection) will be denoted by the same
symbol $m(\wt{\delta})$.

\begin{proposition} \label{pnopieces} Let $\a:[0,1]\to\cc$
be a \crpa path, $\a[0,\tau']$ be its first edge.
For any $\tau\in[0,1)$ let $\wt\delta(\tau)$ be the canonical
representative of the cycle vanishing along the path
$\a|_{[\tau,1]}$. Then
\begin{equation}\log_2m(\wt\delta(0))-\log_2m(\wt\delta(\tau'))
\leq 23n^{12}.\end{equation}
\end{proposition}

The Proposition is proved below.

\begin{corollary} \label{nopieces}  In the previous Proposition
$$m(\wt\delta(0))\leq 2^{E(\a)},\ E(\a)=23n^{12}m(\a), \ m(\a)=\#(\text{edges of}\ \a).$$
\end{corollary}

\begin{proof} {\bf of Proposition \ref{pnopieces}.}
The number of edges of the projection of a canonical representative
increases exactly while crossing a parameter value
$\tau''\in[0,\tau')$
(moving $\tau$ from $\tau'$ towards 0) where
a certain $x_j$ meets some edge of the projection
(the latter edge may be multiple, see the previous Remark). If only
one $x_j$ meets one edge, then this edge
breaks in two pieces, their multiplicities are equal to that of
the initial edge. Therefore,
the number of edges with multiplicities is at most doubled after
passing the value $\tau''$. Analogously,
if there are $r>1$ triples $x_j$, $x_i$,
$x_k$ such that $x_j$ meets the edge $[x_i,x_k]$ at
$\tau=\tau''$, then the total
number of edges with multiplicities increases in at most
$2^r$ times. Hence, the number of edges with multiplicities of the
global vanishing cycle is no greater than $2^N$, where
$N$ is the total number of
quadruples $(x_1,x_2,x_3,\tau)$, $\tau\in[0,\tau')$, such that
$x_i=x_i(\a(\tau))$ form a collinear point triple, and
they are not identically collinear in $\tau\in[0,\tau']$.

Without loss of generality we assume that $\a(\tau)=\tau$
for $\tau\in[0,\tau']$: one can achieve this by
complex affine change of the coordinate $t$. Then the
previous quadruples are isolated solutions
with $\tau\in[0,\tau')$ of the following system of equations

\def\dd{\partial}
\def\re{\operatorname{Re}}
\def\im{\operatorname{Im}}

\begin{equation}
\begin{cases}
H(x_i,y_i)-\tau=0 \\
\frac{\dd H}{\dd y}(x_i,y_i)=0 \\
i=1,2,3 \\
\im\tau=0 \\
(\re x_1-\re x_2)(\im y_1-\im y_3)-
(\re x_1-\re x_3)(\im y_1-\im y_2)=0
\label{hij}\end{cases}
\end{equation}

We show that the number of these
isolated solutions is less than
$23n^{12}$ by using Bezout's theorem applied to
the complexitication of (\ref{hij}), see the following
Definitions.

\begin{definition} \label{defcompl}
The {\it real form} of a system of
polynomials in complex variables is the system of their
real and complex parts
(as polynomials in the real and complex parts of the variables).
\end{definition}
\begin{definition} The {\it complexification} of a system of
real polynomials is the system of their extensions to the complex
variables. The complexification of the real form of a system
(N) will be briefly referred to, as the
complexification of (N).
\end{definition}

The real form of system (\ref{hij}) is
a system of real polynomial equations in $\mathbb R^{14}$
obtained by replacing
the (first 6) complex polynomials by their complexifications.
The complexification of system (\ref{hij}) is a system of
14 polynomial equations in $\cc^{14}$: 6 equations of degree
$n+1$, 6 ones of degree $n$, one linear equation, one quadratic
equation. Bezout's theorem applied to the complexification
says that the number of its isolated solutions is no greater than
the product of the latter degrees
$2(n+1)^6n^6<23n^{12}$ (the latter inequality follows
from elementary inequalities).

\begin{proposition} \label{isola}
Each isolated solution of (\ref{hij}) with
$\tau\in[0,\tau')$ is an isolated solution of its
complexification.
\end{proposition}

\begin{remark} One can provide examples of real
polynomial equations with isolated solutions in the real space
that are not isolated solutions in the complex space:
the real polynomial equation  $x^2+y^2=0$ has unique real
solution 0 that is not an isolated solution of its
complexification. V.Kharlamov have proposed the following example
of 3 real polynomial equations with 3 variables, where
the number of isolated solutions in the real space is
greater than the product of the degrees
(which we call the Bezout number):
$$\begin{cases}
\prod_{k=1}^d(x-x_k)^2+\prod_{k=1}^d(y-y_k)^2=0 \\
z=0 \\
z=0
\end{cases}, 
\ d>2.$$
It has $d^2$ real solutions $(x_i,y_j,0)$, which is greater
than its Bezout number $2d$.
\end{remark}

Proposition \ref{isola} (which is proved below)
together with the previous discussion
implies Proposition \ref{pnopieces}.
\end{proof}

\begin{proof} {\bf of Proposition \ref{isola}.} Fix an isolated
solution $X$ of (\ref{hij}). Let us show that it is an isolated
solution of its complexification. To do this, we
consider the subsystem (denoted
(\ref{hij})') of system (\ref{hij}) consisting of its (first 6)
complex polynomial equations, which is obtained by dropping
its nonholomorphic (two last) equations. Below we show
(Proposition \ref{maxrank}) that the local solutions at $X$ of
system (\ref{hij})' form a regularly embedded
holomorphic curve (two-dimensional
real surface, denoted $\Gamma$) locally
1-to-1 projected to a domain in the
$\tau$- plane. We also show that system (\ref{hij})' has
the maximal rank at $X$. This implies the same statements
in the complexification. In particular, the local solutions at $X$
of the complexified system (\ref{hij})' in $\cc^{14}$ form
a two-dimensional holomorphic surface (denoted $\wt{\Gamma}$)
locally 1-to-1 projected to a
domain in the complexified $\tau$- plane $\cc^2$. Afterwards
we consider the restrictions to $\wt{\Gamma}$ of the complexified
two last equations of (\ref{hij}): it suffices to show
that their solution $X$ is isolated in $\wt{\Gamma}$.
The hyperplane $\operatorname{Im}\tau=0$ is transversal to the real
surface $\Gamma$ at $X$, hence, their intersection
(denoted $\gamma$) is a regularly embedded real analytic curve.
Hence, the two latter statements hold true in the
complexification (denote
$\wt{\gamma}$ the intersection of the complexifications of these
hyperplane and surface, $\wt{\gamma}$
is a regularly embedded holomorphic curve). Now it suffices to show
that the last equation of (\ref{hij}) does not hold identically
on $\wt{\gamma}$. Indeed, it does not hold identically
on the real curve $\gamma$ by isolatedness of
$X$ as a solution of (\ref{hij}). Hence, $X$ is an isolated
solution of the complexification of (\ref{hij}).

Thus, the previous discussion proves Proposition \ref{isola} modulo
the following

\begin{proposition} \label{maxrank}
The solutions of (\ref{hij})' with
$\tau\in [0,\tau')$ form a disjoint union of graphs of analytic
vector functions $(x_i(\tau),y_i(\tau))_{i=1,2,3}$ in
$\tau\in[0,\tau')$.
System (\ref{hij})' has the maximal rank at these solutions.
\end{proposition}

\begin{proof} The solutions
of system (\ref{hij})' with $\tau$ being a noncritical value
of $H$ are exactly the tuples $((x_i,y_i)_{i=1,2,3},\tau)$ where
$(x_i,y_i)$ are the critical points of the projection
$\pi_{\tau}:S_{\tau}\to Ox$.
By assumptions, the path $\a$ is
critically-regular and contains the semiinterval $[0,\tau')$. The
critical regularity implies that the values $\tau\in[0,\tau')$ are
noncritical for $H$, the critical values $x_i(\tau)$
of the projection $\pi_{\tau}$ do not collide
while $\tau$ ranges over the path $\a$ (in particular, when
$\tau\in[0,\tau')$) and are holomorphic functions in $\tau$.
The two last statements imply the same statements with
$x_i(\tau)$ replaced by the critical points
$(x_i(\tau),y_i(\tau))$ of the same projection. This proves the
first statement of the Proposition, which also implies
the transversality of the curves $S_{\tau}\subset\cc^2$
(with $\tau\in[0,\tau')$) and
$$\Gamma_y=\{\frac{\partial H}{\partial y}(x,y)=0\}\subset\cc^2$$
(in particular, the regularity of the latter curve at their
intersection points)

Let us prove the last statement of the Proposition
saying that system (\ref{hij})' has the
maximal rank at $X$. To do this, it suffices to show that
the pair of polynomials $H(x,y)-\tau$
(with fixed $\tau$ and variable $(x,y)$) and
$\frac{\d H}{\d y}(x,y)$ has
the maximal rank at the previous points of intersection of their
zero level curves. By transversality and regularity of the latters,
it suffices to show that the gradients of the previous polynomials
are nonzero at these points. Since
any analytic function with a regular simple zero hypersurface
has a nonzero gradient there, it suffices to show that both
zero level curves are simple. The zero level curve
$S_{\tau}$ of the polynomial $H-\tau$ is simple, since it is
noncritical. The same statement for the zero level curve
$\Gamma_y$ of the polynomial $\frac{\d H}{\d y}$ is implied
by the following

\begin{proposition} \label{simplezeros}
A partial derivative of a ultra-Morse polynomial
has no multiple zero curves.
\end{proposition}

\begin{proof} Let $H$ be a ultra-Morse polynomial.
Without loss of generality we prove
the statement of the Proposition for the derivative
$\frac{\partial H}{\partial y}$. Suppose the contrary:
there exists an irreducible algebraic curve $\Gamma'_y$
where the derivative has multiple zero. Let us show that then
$H$ has a critical point that is not Morse. The contradiction
thus obtained will prove the Proposition.

Consider the other
partial derivative $\frac{\partial H}{\d x}$ restricted to
$\Gamma_y'$. We claim that $\frac{\d H}{\d x}(Z)\to\infty$, as
$Z\to\infty$ along $\Gamma_y'$. This follows from the fact that
the intersection of $\Gamma_y'$ with the infinity line (in the
projective plane) is disjoint from the intersection of this line
with zero level curve of $\frac{\d H}{\d x}$.
Indeed, the latter intersections remain the same, if $H$
is replaced by its highest homogeneous part $h$. The
 partial derivatives of $h$ are relatively prime, thus,
have no common zero lines ($h$ is generic, since $H$ is
ultra-Morse).

Thus, $\frac{\d H}{\d x}$ is a holomorphic function on the curve
$\Gamma_y'$ with pole(s) at infinity. Therefore, it has at least
one zero on this curve. This zero is a critical point of $H$ by
definition. It is not Morse, since $\frac{\d H}{\d y}$ has
zero gradient on the ambient curve $\Gamma'_y$ (where it has
multiple zero). This proves the Proposition.
\end{proof}

Proposition \ref{simplezeros} together with the previous discussion
implies that system (\ref{hij})' has the maximal rank at $X$.
This proves Propositions \ref{maxrank} and \ref{isola}.
\end{proof}
\end{proof}

\subsection{Lengths of canonical representatives of vanishing
cycles}

\begin{lemma} \label{l(n)}
Let $H$ be an ultra-Morse polynomial of
degree $n+1\geq3$, $a$ be its critical value,
$t\in\mathbb C$ be a noncritical value, $S_t=\{H(x,y)=t\}$,
$\pi_t:S_t\to Ox$
be the projection to the $x$- axis. Let $\alpha$ be a
critically-regular path from $t$ to $a$ with $m(\alpha)$
edges, $\wt\delta\subset S_t$ be the canonical
representative of the corresponding vanishing cycle.
Let $X,Y>0$ be such that the topology of the curve $S_t$
be contained in the bidisc $D_{X,Y}$. Then
\begin{equation}|\wt\delta|< 2^{l(n)m(\alpha)}R,\ R=\max\{ X,Y\},\
l(n)=24n^{12}.\label{rln}\end{equation}
\end{lemma}

\begin{proof}
Each arc couple of $\wt\delta$ (see   Proposition \ref{pieces}),
which is  projected onto a real segment in the complex $Ox$ -line,
lies in a real algebraic curve. This curve (denoted
$\psi$) is the intersection of
the complex level curve $S_t$ of $H$  and the real
hyperplane in $\cc^2=\mathbb R^4$
that contains the complex $Oy$- line and the latter segment.
That is, $\psi$ is the common zero set of the system
of 3 polynomials in $\cc^2=\mathbb R^4$: the real and the imaginary
parts of $H$ (both of degree $n+1$) and a linear function.

The number of the previous arc couples (edges)
is already estimated in Corollary
\ref{nopieces}. Let us estimate the total length
of one arc couple. It is no greater than the
sum of lengths of its projections to the $x$- axis and to the
real and imaginary $y$ -axes (more precisely, we have to replace
"length of projection" by "length of projection times the maximal
number of preimages of a generic point").
Let us estimate the latter quantities.

The arc couple is projected to the $x$- axis onto
a segment in $D_X$ (hence, having length less than $2X$), and each
point of the latter segment (except for its ends)
has exactly two  preimages (Proposition \ref{pieces}).
Hence, the length contribution of the $x$- projection is
less than $4X$.

The projection image of the arc couple to either real or imaginary
$y$- axis is a segment lying in $D_Y$ (hence, of length less than
$2Y$). The number of preimages of a generic point
is no greater than $(n+1)^2$. Indeed, it suffices to show that
the ambient algebraic curve $\psi$
has at most the same number of isolated intersection points
with a generic real hyperplane in $\cc^2=\mathbb R^4$
parallel to a given one (here "generic" means "intersecting
$\psi$ transversally"). These points are common zeros
of the previous system of 3 polynomials (defining $\psi$) and
an additional linear function (defining the hyperplane).
 By transversality, these points are
isolated common zeros of their complexifications. Hence,
by Bezout's theorem, their number is no  greater than $(n+1)^2$.

Therefore, the length contribution of
the projections to either real or imaginary $y$- axis is less than
$2Y(n+1)^2$.

The previous discussion implies that the total length of an arc
couple is less than
$$4X+4Y(n+1)^2\leq (4+9n^2)Y\leq10Rn^2.$$
Together with Corollary \ref{nopieces}, this implies that
the total length of the canonical representative is less than
$$10Rn^2\times 2^{23n^{12}m(\a)}<2^{24n^{12}}m(\a)R.$$
This proves Lemma \ref{l(n)}.
\end{proof}

We use also the following Corollary of Lemma \ref{l(n)} giving an
upper bound of length
of vanishing cycle in terms of the length of the corresponding path.
To formulate it, let us introduce the next Definition.

\begin{definition} \label{dbreg}
Let $H$ be a ultra-Morse polynomial, $\beta>0$.
A regular path $\a:[0,1]\to\cc$ (denote $a=\a(1)$) is said to be
$\beta$- {\it regular}, if the curve $\a\cap D_{\beta}(a)$ is
a connected arc of the path $\a$,
and $\a$ is disjoint from the $\beta$ neighborhoods of the
critical values distinct from $a$ of the polynomial $H$.
\end{definition}

\begin{corollary} \label{cl(n)}
Let $H$ be an ultra-Morse polynomial of degree
$n+1\geq3$, $0<\beta<1$. Let  $a$ be its critical value,
$t\in\mathbb C$ be a noncritical value,
$\alpha$ be a $\beta$- regular (but now, not necessarily
piecewise-linear) path from $t$ to $a$.
 Let $X,Y>0$ be such that the topology of the curve $S_t$
lies in the bidisc $D_{X,Y}$, $R=\max\{ X,Y\}$. Then
the cycle in $H_1(S_t,\mathbb Z)$ vanishing along the path
$\alpha$ admits a representative $\wt{\delta}$ such that
\begin{equation}
\wt{\delta}\subset D_{X,Y},\ |\wt{\delta}|<
2^{l(n)\frac{|\alpha|\beta^{-1}+5}3}R,\ \text{where}
\label{lengthup}\end{equation}
$l(n)$ is the same, as in (\ref{rln}).
\end{corollary}

\begin{proof} The vanishing cycle depends only on the homotopy class
(modulo the critical values) of the path $\alpha$. We construct
a critically-regular path $\a'$ homotopic to $\a$ with at most
$\frac{|\a|\beta^{-1}+5}3$ edges . Then Lemma \ref{l(n)} applied to
$\a'$ implies (\ref{lengthup}) for the canonical representative
$\wt{\delta}$ of the cycle as vanishing along $\a'$. To do this,
consider a splitting of $\alpha$ into its arc $\a\cap D_{\beta}(a)$
(whose length is no less than $\beta$) and
at most $\frac{|\alpha|\beta^{-1}+2}3$ other arcs
of lengths at most $3\beta$. Each splitting arc is homotopic outside the
critical values of $H$ to the straightline segment with the
same ends. Indeed, by $\beta$- regularity,
the $\beta$- neighborhood of the path end $a$ does not contain
other critical values of $H$. This implies that the
arc $\a\cap D_{\beta}(a)$ is homotopic to the radius with the same ends.
The similar statement for the other arcs, which are disjoint
from the $\beta$- neighborhoods of all the critical values,
is implied now by the following
geometric fact: a curve lying outside a disc (or several discs) of
a radius $\beta$ and having a length less than $\pi \beta>3\beta$
is always homotopic  with fixed ends (outside the centers of the discs)
to the straightline segment with the same ends. Therefore, the
new path (denoted $\alpha'$) from $t$ to $a$ consisting of
the previous straightline segments is homotopic to $\alpha$ (the
homotopy does not
meet the critical values of $H$). The path $\alpha'$
is piecewise-linear with at most $\frac{|\alpha|\beta^{-1}+5}3$
edges. One can make it critically-regular by small deformation
of its vertices. Then this is a path $\a'$ we are looking for.
This together with the previous discussions proves the Corollary.
\end{proof}

We will also use the following
upper bound of length of vanishing cycle along
a path that is not $\beta$- regular, in particular,
intersects the $\beta$- neighborhoods of
the critical values. The bound is given in terms of
the length of its part that lies in the disc
$\overline D_3$ outside these neighborhoods and
the variations of arguments of $t-a_i$ (respectively, $t$)
along its arcs in $D_{\beta}(a_i)$ (respectively, outside
$\overline D_3$).

\def\Var{\operatorname{Var}}

\begin{corollary} \label{c2ln} Let $H$ be a normalized ultra-Morse
polynomial of degree $n+1\geq3$, $a_i$ be its critical values.
 Let $0<\beta\leq\nu=\frac{c''(H)}{4n^2}$, $\a:[0,1]\to\cc$
be a regular path,
$$t_0=\a(0),\ a=\a(1),\ \tau'=\min\{\tau\in[0,1]\ | \ \a(\tau,1]
\subset D_{\beta}(a)\}, \ \hat\a=\a\setminus\a(\tau',1],$$
\begin{equation}\mathcal D^{\beta}=
\overline D_3\setminus\cup_iD_{\beta}(a_i), \
\wt\a=\a\cap\mathcal D^{\beta},\label{cald}\end{equation}
$$
V=V_{\a,\beta}=\beta\sum_i\Var_{\hat\a\cap D_{\beta}(a_i)}\arg(t-a_i)
+3\Var_{\hat\a\setminus\overline D_3}\arg t$$
(the $\Var$- terms are the complete variations of arguments along
the corresponding pieces of $\a$). Let $X,Y>0$ be such
that the topology of the curve $S_{t_0}$
lies in the bidisc $D_{X,Y}$, $R=\max\{ X,Y\}$. Then
the cycle in $H_1(S_{t_0},\mathbb Z)$ vanishing along the path
$\alpha$ admits a representative $\wt{\delta}$ such that
\begin{equation}
\wt{\delta}\subset D_{X,Y},\ |\wt{\delta}|<
2^{2l(n)\frac{(|\wt\a|+V)\beta^{-1}+5}3}R,\ l(n)=24n^{12}.
\label{lengthup2}\end{equation}
\end{corollary}

\begin{proof} Recall that each disc
$D_{\beta}(a_i)$ contains no critical values of $H$
except for $a_i$, and these discs are disjoint,
since $\beta\leq\nu$. The complement
$\cc\setminus D_3$ also contains no critical values and is disjoint
from the latter discs by normalizedness  ($|a_i|\leq2$).

If $\hat\a=\emptyset$, then $\a\subset D_{\beta}(a)$, and hence,
$\a$ is homotopic to the segment $[t_0,a]$ modulo the critical
values.  In this case Lemma
\ref{l(n)} applied to the segment path $[t_0,a]$
implies an upper bound stronger than (\ref{lengthup2}).

If $\a$ is contained in $\overline D_3$ and is $\beta$- regular,
then Corollary \ref{cl(n)} implies
an upper bound stronger than (\ref{lengthup2}).

In what follows we assume that $\a$ is either
not $\beta$- regular, or not contained in
$\overline D_3$.
We construct a $\beta$- regular path $\a'':[0,1]\to\overline D_3$
of length at most $|\wt\a|+V$ such that the composition
$\phi=[t_0,\a''(0)]\circ\a''$ is homotopic to $\a$ modulo the
critical values. Then we replace the path $\a''$ by
a homotopic  piecewise linear path $\a'$ with at most
$\frac{|\a''|\beta^{-1}+5}3$ edges, as in the proof of
Corollary \ref{cl(n)}. The path $\phi'=[t_0,\a''(0)]\circ\a'$
thus obtained, which is
homotopic to $\a$, is a piecewise linear path with one edge more.
Applying Lemma \ref{l(n)} to $\phi'$ yields (\ref{lengthup2}).

The previous path $\phi$ is constructed by replacing
connected components of the complement
$\a\setminus\mathcal D^{\beta}$.
Let $l=\a(\tau_1,\tau_2)$ be a maximal arc of
$\a$ lying in one of the
latter connected components (maximal means that the previous
interval $(\tau_1,\tau_2)$ is contained in no other interval
with the same property): then either
$l\subset D_{\beta}(a_i)$ for some $i$, or
$l\subset\cc\setminus\overline D_3$.
If $l=\a(\tau',1]$, we do not replace it. If
not, then $l$ is either
a starting arc $\a[0,\tau'')$ of the path $\a$ with one end
$\a(\tau'')$ in $\partial\mathcal D^{\beta}$, or an arc with
both ends in one and the same boundary circle. In the latter
case we replace $l$ by the arc (denoted $l'$)
of this circle with the same ends
that is homotopic to $l$ outside the critical values
(it may be self-overlapped). Then it
follows from definition that
$$|l'|\leq\beta\Var_{l}\arg(t-a_i),\ \text{if}\ l\subset
D_{\beta}(a_i);\ |l'|\leq3\Var_l\arg t,\ \text{if}\ l\cap D_3=
\emptyset.$$

Now consider the former case, when $l$ is a starting arc
(either inside some $D_{\beta}(a_i)$, or outside $D_3$).
Then we replace $l$ by the composition of

- the segment joining $\a(0)$ to
the closest point of the corresponding circle
(either $\partial D_\beta(a_i)$, or $\partial D_3$);

- an  arc of the latter circle

so that their composition be  homotopic to $l$ in the punctured
disc $D_{\beta}(a_i)$ (respectively, in the complement of $D_3$).

Then it follows from construction that the length of the latter
circle arc is no greater than $\beta\Var_l\arg(t-a_i)$
(respectively, $3\Var_l\arg t$).

The previous replacements give us a modified path $\a$
(denoted $\phi$) that either starts in $\mathcal D^{\beta}$
(then we put $\a''=\phi$), or starts outside by a segment
ending at $\partial\mathcal D^{\beta}$ (then we put
$\a''$ to be $\phi$ with the latter segment deleted).
It follows from construction  that $\a''$ is $\beta$- regular.
The previous inequalities imply that
$|\a''|\leq|\wt\a|+V$. The Corollary is proved.
\end{proof}

\subsection{Intersection indices of vanishing cycles}
\begin{theorem} \label{tindex} Let $H$ be an ultra-Morse polynomial,
$t\in\cc$ be a noncritical value, $\a_1$, $\a_2$ be two
critically-regular paths starting at $t$ and going to some
critical values of $H$ (may be coinciding). Then the
module of the intersection
index of the corresponding vanishing cycles is less than
$$2^{24n^{12}(m(\a_1)+m(\a_2))}$$
\end{theorem}

\begin{proof} Let $\wt\delta_1$, $\wt\delta_2$ be the canonical representatives
of  the vanishing cycles. Their projections
are piecewise-linear curves with number of edges estimated
from above by Corollary \ref{nopieces}.
For the proof of the Theorem we estimate
the number of the intersection points of the projections and
the contribution of each point to the intersection index.

Case 1. Two arc couples of $\wt\delta_1$ and $\wt\delta_2$
respectively have transversally intersected projection edges
(no common vertex). The contribution of this intersection
point (per arc couple pair) to the
intersection index of $\wt\delta_1$ and $\wt\delta_2$
has module at most 2. Indeed, each arc couple
consists of two arcs (disjoint maybe except for their ends),
each arc is 1-to-1 projected to
the corresponding edge (Proposition \ref{pieces}). Hence, we have
at most two transversal intersection points
of the corresponding unions of arcs over the intersection point.

There can be at most $2^{E(\a_1)+E(\a_2)}$ arc couple pairs
(one arc couple from $\wt\delta_1$, the other one from $\wt\delta_2$)
with transversally intersected projections, as above.
This follows from Corollary \ref{nopieces}.
Thus, the total contribution of the transversal intersections
of edges to the intersection index
is no greater than $2^{E(\a_1)+E(\a_2)+1}$.

Case 2. There is a common vertex $x$ of the projections
$\pi_t\wt{\delta_1}$ and $\pi_t\wt{\delta_2}$ (with at most 2 adjacent
edges in each $\pi_t\wt{\delta_i}$; thus, the total number
of adjacent edges of both projections is at most 4).
Firstly let us assume the latter edges are simple and no two of them
are collinear. We claim that the module of the contribution of $x$ to the
intersection index of the canonical representatives
is at most 16. Indeed, by assumption, each projection $\pi_t\wt{\delta_i}$ has two
adjacent edges at $x$, and $\wt{\delta_i}$ has two arcs over each edge.
Each arc is locally (in a neighborhood of the preimage
$\pi_t^{-1}(x)$) a regularly embedded semicurve with a tangent line
at its end over $x$ (its end may be either a
critical point of projection or not); in total there are
4 semicurves per each
$\pi_t\wt{\delta_i}$. All these 8 semicurves are transversal by
construction. This implies that the index of the
intersection over $x$ of the two unions of 4 semicurves is no
greater than $4\times4=16$.

Let now there be a common vertex of projections (with at most two
adjacent edges in each $\pi_t\wt{\delta}_i$, as above) and either
some of the (at most four) adjacent edges are collinear,
or some edges coincide. Then one can avoid any of these situations
by small deformation of the interiors of edges (keeping
the vertices fixed) in class of smooth curves and  lifting it
up to a deformation of the closed curves $\wt{\delta}_i$.
If the deformation of edges is $C^1$- small enough, then the number
of transversal intersection points of projections (not including
vertices) remains the same.  Afterwards the contribution
to the intersection index of each common vertex
of the projections is estimated as above: it is no greater than 16.

In the general case, when there are several adjacent edges
and (or) some of them are multiple,
define the multiplicity of a vertex of the projection
$\pi\wt{\delta_i}$ as
the total multiplicity of the adjacent edges of the same projection.
The total vertex multiplicity of the projection is the sum of the
multiplicities of its vertices.
It follows from definition and the previous discussion
that if $A$ is a common vertex of
$\pi\wt{\delta}_i$, $i=1,2,$ with corresponding multiplicities $\mu_i$,
then its contribution to the intersection index is no greater than
$16\mu_1\mu_2$. Hence, the total contribution of the common
projection vertices to the intersection index
is no greater than 16 times the product of the total vertex multiplicities of the  projections.

Each edge has two vertices, hence, the total multiplicity of vertices
in each
$\pi\wt{\delta}_i$ is no greater than the double total multiplicity
of edges, thus, no greater than $2^{E(\a_i)+1}$
(Corollary \ref{nopieces}).
Therefore, the total  contribution of the common
vertices to the intersection index is no greater than
$2^{E(\a_1)+E(\a_2)+6}$.

The previous discussion implies that the intersection index of the
vanishing cycles under consideration, which is the sum of
contributions of transversal intersections of the projections
and those of their common vertices, is no greater than
$2^{E(\a_1)+E(\a_2)+7}$.
Together with Corollary \ref{nopieces} and elementary
inequalities, this implies that the module of the
intersection index is less than $2^{24n^{12}(m(\a_1)+m(\a_2))}$.
Theorem \ref{tindex} is proved.
\end{proof}

\subsection{Upper bounds of integrals}
Below we recall and prove upper bounds of Abelian integrals
used in \cite{gi}.
To state them, let us firstly recall the definition (introduced
in \cite{gi}) of the following function
on the space of normalized ultra-Morse polynomials.

\begin{definition} Let $H$ be a normalized ultra-Morse polynomial
of degree $n+1\geq3$.
Define $c_2(H)$ to be $n^2$ times the minimal distance between
its critical values. Put
$$c''(H)=\min(c_2(H),1).$$
\end{definition}

\begin{theorem} \label{t2.10} Let $H$ be a normalized ultra-Morse
polynomial of degree $n+1\geq3$,
$a$ be its critical value, $t$ be a noncritical value,
$|t|\leq5$, $\a:[0,1]\to\cc$ be a $\varepsilon$-
regular path from $t$ to $a$ (See Definition \ref{dbreg}),
\begin{equation} |\var|=\frac{c''(H)}{8n^2},\ |\a|\leq36n^2+10.
\label{|a|<}
\end{equation}
Let $\delta\in H_1(S_t,\zz)$ be the
cycle vanishing along $\a$. Let $\omega$ be a monomial 1- form
of degree no greater than $2n-1$ with unit coefficient. Then
\begin{equation} |I_{\delta}(t)|=|\int_{\delta}\omega|
<2^{\frac{2600n^{16}}{c''(H)}}(c'(H))^{-28n^4}.
\label{uppi1}\end{equation}
\end{theorem}

\begin{theorem} \label{t3.1} Let $H$ be a weakly normalized ultra-Morse
polynomial of degree $n+1\geq3$, $a$ be its critical value, $t$ be a noncritical value, $0<\beta<1$.
Let $\a$ be a $\beta$- regular path from $t$ to $a$ (see Definition
\ref{dbreg}),
$\delta\in H_1(S_t,\zz)$ be the corresponding vanishing cycle.
 Let $\omega$ be a monomial
1- form of degree at most $2n-1$ with unit coefficient. Then
\begin{equation} |I_{\delta}(t)|<2^{10n^{12}\frac{|\a|+5}{\beta}-
2n}
(c'(H))^{-28n^4}.\label{uppi2}\end{equation}
\end{theorem}

\begin{theorem} \label{tpoinc} Let $H$ be a normalized ultra-Morse
polynomial of degree $n+1\geq3$, $a$ be its critical value,
$t_0$ be a noncritical value, $0<\beta\leq\nu=\frac{c''(H)}{4n^2}$.
Let $\a:[0,1]\to\cc$ be a regular path, $t_0=\a(0)$, $a=\a(1)$,
$\hat\a$, $\wt\a$, $V=V_{\a,\beta}$ be the same, as in (\ref{cald}).
Let $\delta\in H_1(S_{t_0},\zz)$ be the cycle vanishing along $\a$.
 Let $\omega$ be a monomial
1- form of degree at most $2n-1$ with unit coefficient. Then
\begin{equation} |I_{\delta}(t_0)|<
2^{20n^{12}\frac{|\wt\a|+V+5}{\beta}-2n}
(c'(H))^{-28n^4}\max\{1,(\frac{|t_0|}5)^2\}.\label{uppi3}\end{equation}
\end{theorem}

\begin{proof} {\bf of Theorem \ref{t2.10}.}
Let us apply Theorem \ref{topxy} to the polynomial $H$.

Case 1: the coordinates $(x,y)$ satisfy the statements of Theorem
\ref{topxy}. Then this Theorem implies that the topology of
the curve $S_t$ lies in the bidisc $D_{X,Y}$, $X\leq Y\leq R_0$,
$R_0$ is the same, as in (\ref{xy<<}). By assumption, the
conditions of Corollary \ref{cl(n)} hold with $\beta=\var$,
$R\leq R_0$. Let $\wt{\delta}\subset S_t\cap D_{X,Y}$ be the
representative from this Corollary of the vanishing cycle. Then
by the same Corollary  and definition,
\begin{equation} |\int_{\wt{\delta}}\omega|\leq R_0^{2n-1}
|\wt{\delta}|\leq
R_0^{2n}2^{l(n)(\frac{|\a|+5\var}{3\var})},
\ l(n)=24n^{12}.\label{iwd1}\end{equation}
Substituting the values $\var=\frac{c''(H)}{8n^2}$,
$R_0=(c'(H))^{-14n^3}n^{65n^3}$
and  inequality (\ref{|a|<}) to the latter right-hand side
and applying elementary inequalities yields
\begin{equation}|\int_{\wt{\delta}}\omega|\leq
(c'(H))^{-28n^4}2^{\frac{2600n^{16}}{c''(H)}-4n^{16}}.
\label{iwd2}\end{equation}
This proves Theorem \ref{t2.10} in Case 1.

\medskip

Case 2: general. Let $(x',y')$ be orthogonal coordinates on $\ccd$
satisfying the statements of Theorem \ref{topxy}. The form $\omega$
is monomial of degree at most $2n-1$ with unit coefficient.
Therefore, in the new coordinates it becomes a product
of a constant 1- form $A_1dx'+B_1dy'$ and at most
$2n-1$ linear functions $A_ix'+B_iy'$, $i\geq2$; the latter
form and linear functions have unit Hermitian norm:
$|A_i|^2+|B_i|^2=1$.

The sum of modules of the coefficients of
the form $\omega$ in the new coordinates
is no greater than $2^n$. Indeed, it is no greater than
$$\prod_i(|A_i|+|B_i|),\ |A_i|+|B_i|\leq2\sqrt{\frac{|A_i|^2+
|B_i|^2}2}=\sqrt2$$
(the classical quadratic mean inequality). Hence,
the previous product
(and thus, the sum of modules of coefficients) is
no greater than $2^n$.

Let us repeat the discussion from Case 1 in the coordinates
$(x',y')$. Now our form is not necessarily monomial, and
the previous upper bounds of the integral should be multiplied
by the sum of modules of its coefficients. This together
with the previous inequality implies the same upper bound
(\ref{iwd2})  but with the right-hand side multiplied by $2^n$.
The right-hand side thus modified is less than that of the
inequality in Theorem \ref{t2.10}. This proves Theorem \ref{t2.10}.
\end{proof}

\begin{proof} {\bf of Theorem \ref{t3.1}.}  By assumptions,
all the critical values of $H$ have modules at most 2. Hence,
$$|t|\leq|\a|+2.$$
Let us assume that $|t|>5$, thus, $|\a|+2>5$
(the opposite case is treated
simpler and we get a stronger inequality than in the Theorem;
this case will be briefly discussed at the end of the proof).
We consider only the case when the coordinates in $\ccd$ satisfy
the statements of Theorem \ref{topxy}: then afterwards
one has only to check that the upper
bound obtained remains less than that in Theorem \ref{t3.1}
after multiplication by $2^n$. This will imply the Theorem
in the general case, as in the proof of Theorem \ref{t2.10}.

By the Corollary of Theorem \ref{topxy} and the previous
inequality on $|t|$, the topology of the
curve $S_t$ is contained in a bidisc $D_{X,Y}$,
\begin{equation}X\leq Y\leq R=R_0(\frac{|\a|+2}5)^{\frac1{n+1}},
\ R_0\ \text{is the same as in Theorem \ref{topxy}.}
\label{a+2}\end{equation}
Then as in (\ref{iwd1}), by definition and Corollary \ref{cl(n)},
$$|\int_{\delta}\omega|\leq R^{2n}2^{l(n)\frac{|\a|+5\beta}{3\beta}}, \
l(n)=24n^{12}. \ \text{Thus,}$$
$$|\int_{\delta}\omega|< R_0^{2n}(\frac{|\a|+2}5)^2
2^{24n^{12}\frac{|\a|+5\beta}{3\beta}}.$$
Substituting the value of $R_0$ from Theorem \ref{topxy} together
with elementary inequalities yields
\begin{equation}|\int_{\delta}\omega|<2^{10n^{12}
\frac{|\a|+5}{\beta}-3n}
(c'(H))^{-28n^4}(\frac25)^2.\label{a+2n}\end{equation}
The latter right-hand side (even being multiplied by $2^n$) is
less than that of (\ref{uppi2}). This proves (\ref{uppi2}).

Let us now consider the case, when $|t|\leq5$. Then the
previous upper bounds of the integral hold, but now we have
to replace the ratio $\frac{|\a|+2}5$ in (\ref{a+2})
and the next inequality by 1. This implies that (\ref{a+2n})
holds with
its right-hand side multiplied by $(\frac5{|\a|+2})^2<(\frac52)^2$.
Thus modified right-hand side (even if multiplied by $2^n$)
is again no greater than
that of (\ref{uppi2}). Theorem \ref{t3.1} is proved.
\end{proof}

\begin{proof} {\bf of Theorem \ref{tpoinc}.}
The proof of Theorem \ref{tpoinc} repeats that of Theorem \ref{t3.1}
with obvious change: we have to substitute the length estimate
of vanishing cycle given by Corollary \ref{c2ln}
(instead of \ref{cl(n)})
and $R=R_0\max\{1,(\frac{|t_0|}5)^{\frac1{n+1}}\}$ (by Corollary
\ref{topxyt}).
\end{proof}

\subsection{Lower bound of period determinant}

\begin{definition} (see \cite{gi}) A tuple $\Omega=(\o_1,\dots,\o_{n^2})$
of monomial 1- forms $\o_i$ of the type $x^ly^{m+1}dx$
is called {\it standard}, if

-  their degrees are no
greater than $2n-1$.

- all the forms $x^ly^{m+1}dx$, $0\leq l+m\leq n-1$, are contained
there.
\end{definition}

The following Theorem was stated and used in \cite{gi}.

\begin{theorem} \label{lowerdet} Let
$H$ be a normalized ultra-Morse polynomial,
$\Omega$ be a standard monomial tuple of forms.
Let $\mathbb I(t)$, $\D(t)$ be respectively the
corresponding Abelian integral matrix (\ref{2.3})
and its determinant (\ref{dett}).
The standard monomial tuple $\Omega$ can be chosen so  that
\begin{equation}|\D(t)|>(c'(H))^{6n^3}(c''(H))^{n^2}n^{-62n^3},
\label{1.16}
\end{equation}
whenever $t$ lies outside the $\frac{c''(H)}{4n^2}$- neighborhoods
of the critical values of $H$.
\end{theorem}

The proof of this theorem is based on an explicit formula (\ref{2.20}):
$$\D(t)=C(h,\Omega)\prod_{i=1}^{n^2}(t-a_i),$$
 for the period determinant, see \cite{g} and \cite{ga},
which holds, if
the form tuple under consideration is $d$- standard, see the following Definition.

\begin{definition} A standard monomial form tuple is called $d$-
{\it standard}, if for any $d\in\mathbb N$, $n\leq d\leq2n-2$, it
contains exactly $2n-d-1$ forms of degree $d+1$.
\end{definition}

\begin{example} Let $(l(i),m(i))$ be a lexicographic sequence of pairs
$(l,m)$, $0\leq l,m\leq n-1$, numerated by $i=1,\dots,n^2$. Put
$$\omega_i=x^{l(i)}y^{m(i)+1}dx.$$
This is a $d$- standard form tuple, which follows from definition.
If in the previous Theorem the highest homogeneous part of $H$ equals
$h(x,y)=x^{n+1}+y^{n+1}$, then the latter forms satisfy its inequality.
The proof of this statement, which is omitted to save the space, can be easily derived from the results of this Subsection and Section 3.
In general, one can choose a generic $h$ in such a way that
the determinant $\Delta_{h,\Omega}(t)$ constructed with the above
$\omega_i$
be identically equal to 0 (this follows from results of \cite{g}):
in this case the form tuple should be changed.
\end{example}

An explicit formula for the constant $C(h,\Omega)$ was obtained in \cite{g}.
\begin{equation}
C(h,\Omega)=C_n(\Sigma(h))^{\frac12-n}P(h,\Omega),\ \text{where}\label{2.21}
\end{equation}
$\Sigma(h)$ is the discriminant of the homogeneous polynomial $h$:
\begin{equation}\text{if}\ h(x,y)=c_{-1}\prod_{i=0}^{n}(y-c_ix), \ \text{then}\
\Sigma(h)=c_{-1}^{2n}\prod_{0\leq j<i\leq n}(c_i-c_j)^2,\label{2.4}
\end{equation}

$$P(h,\Omega)=\prod_{d=n}^{2n-2}P_d(h,\Omega),$$
where $P_d$ are the polynomials from \cite{g},
whose definition is recalled below,

\begin{equation}
C_n=(-1)^{\frac{n(3n-1)}4}\frac{(2\pi)^{\frac{n(n+1)}2}(n+1)^{
\frac{n^2+n-4}2}((n+1)!)^n}{\prod_{m=1}^{n-1}(m+n+1)!}.\label{2.9}
\end{equation}

\begin{definition}\label{dome} (\cite{g}). For any given
polynomial 1- form $\omega$
put $D\omega$ to be the polynomial defined by the equality
$$d\omega=D\omega dx\wedge dy.$$
\end{definition}

\begin{definition} \label{defado} Let   $n\geq2$, $d\in\mathbb N$,
$n\leq d\leq 2n-2$,
$h$ be a homogeneous polynomial of degree $n+1$. Let $\Omega(d)=(\o_1',\dots,\o_s')$ be
an ordered tuple of homogeneous 1- forms of degree $d+1$,
the number $s$ of the forms being
equal to $s=d+1$ in the case, when $d\leq n-1$, and
$s=2n-d-1$ otherwise. The {\it matrix $A_d(h,\Omega(d))$ associated
to the form tuple $\Omega(d)$} is
the $(d+1)\times(d+1)$ matrix whose columns are numerated   by all the
monomials $y^d, y^{d-1}x,\dots,x^d$ of degree $d$ and the lines consist of
the corresponding coefficients of the following polynomials:

Case $d\leq n-1$. Take the $d+1$ polynomials $\frac{D\omega_r'}{d-r+2}$.

Case $d\geq n$. Take the
$d-n+1$ polynomials $x^jy^{d-n-j}\frac{\partial h}{\partial y}$,\
$0\leq j\leq d-n$; the $2n-d-1$ polynomials $\frac{D\omega_r'}{n-r+1}$;
 the $d-n+1$ polynomials $x^jy^{d-n-j}
\frac{\partial h}{\partial x}$, $0\leq j\leq d-n$.
\end{definition}

Let $\Omega$ be a $d$- standard form tuple,
$\Omega(d)$ be the tuple of the
forms in $\Omega$ of degree $d+1$ (numerated in the same order, as in $\Omega$).
The number $s$ of forms in $\Omega(d)$ is equal to $2n-d-1$ by definition.  Put
\begin{equation}
A_d(h,\Omega)=A_d(h,\Omega(d)),\ P_d(h,\Omega)=
P_d(h,\Omega(d))=\det A_d(h,\Omega(d)),\label{2.19}
\end{equation}
$$P(h,\Omega)=\prod_{d=n}^{2n-2}P_d(h,\Omega).$$
Now all the entries of formula (\ref{2.21}) are defined.

As it will be shown below, Theorem \ref{lowerdet} is implied by the following
\begin{theorem} \label{tlowerch} Let $h$ be a given generic polynomial
of degree $n+1\geq3$ with $||h||_{max}=1$.
Then one can choose a $d$- standard form tuple $\Omega$ so that
\begin{equation}C(h,\Omega)>(c'(h))^{6n^3}n^{-60n^3}.\label{lowerch}
\end{equation}
\end{theorem}

Theorem \ref{tlowerch} will be proved in Section 3.
The principal part of its proof is the lower bound of $P_d(h,\Omega)$ given by the
next Lemma. Its proof takes the most of Section 3.
The upper bound of $\Sigma(h)$ and lower bound of $C_n$ are proved in Section 3 using
elementary inequalities and straightforward a priori estimates for $h$.

\begin{lemma} \label{lpd>} For any generic homogeneous polynomial $h$
of degree $n+1\geq3$ with $||h||_{max}=1$
and any $d=n,\dots,2n-2$ one can choose
a collection $\Omega(d)=(\omega_1,\dots,\omega_s)$ of $s=2n-d-1$
 monomial  forms of the type $x^ly^{m+1}dx$, $l+m=d$, so that
\begin{equation} P_d(h,\Omega(d))> n^{-44n^2}(c'(h))^{6n^2}.
\label{pd>}\end{equation}
\end{lemma}

\begin{proof} {\bf of Theorem \ref{lowerdet}.} Let $t\in\cc$ lie outside
$\frac{c''(H)}{4n^2}$- neighborhoods of the critical values of $H$.
Then by (\ref{lowerch}) and (\ref{2.20}),
$$|\Delta(t)|>(c'(h))^{6n^3}n^{-60n^3}(\frac{c''(H)}{4n^2})^{n^2}.$$
This together with elementary inequalities implies (\ref{1.16}) and proves
Theorem \ref{lowerdet} modulo Theorem \ref{tlowerch}.
\end{proof}

\section{Upper bounds of topology. Proof of Theorems \ref{topxy}
and \ref{t3.5}}

In this Section we give a proof of Theorem \ref{t3.5} (the most part of
the Section). In 2.2 we prove Theorem \ref{topxy} and its Addendum.
In 2.3 we prove Proposition \ref{phl} and some more precise a
priori bounds
for unit-scaled polynomials. These bounds are used in the proof  of Theorems
\ref{topxy} and \ref{t3.5}.

\subsection{The plan of the proof of Theorem \ref{t3.5}}

It suffices to prove Theorem \ref{t3.5} for ultra-Morse polynomials.
The same statement in the general case then follows by passing to non ultra-Morse
limit. Thus, everywhere below (except for Subsection 2.2 and
whenever the contrary is specified)
without loss of generality we consider that the polynomial $H$ is ultra-Morse.

In the proof of Theorem \ref{t3.5} we use the following well-known
and elementary topological properties of basic cycles in $H_1(S_t,\zz)$.
\begin{proposition} Let $H$ be an ultra-Morse polynomial, $t\in\cc$ be a noncritical
value, $\gamma\subset S_t$ be an embedded curve that joins two distinct points
at the infinity line of the compactified curve $\overline S_t\subset\mathbb{CP}^2$.
Then there exists a cycle in the homology group of the affine curve $S_t$ that has
a nonzero intersection index with $\gamma$ (see Fig.2).
\end{proposition}

\begin{proof} Take a cycle close to infinity and surrounding an
end of $\gamma$.
\end{proof}

\begin{corollary} \label{corint} Given a ultra-Morse polynomial $H$, a noncritical value
$t\in\mathbb C$ and a set of generators in $H_1(S_t,\zz)$. Let $b\in Ox$ be a critical
value of the projection $\pi_t:S_t\to Ox$. Then the projection of any representative
of some generator intersects any ray issued from $b$, see Fig.2.
\end{corollary}

\begin{figure}[ht]
  \begin{center}
   \epsfig{file=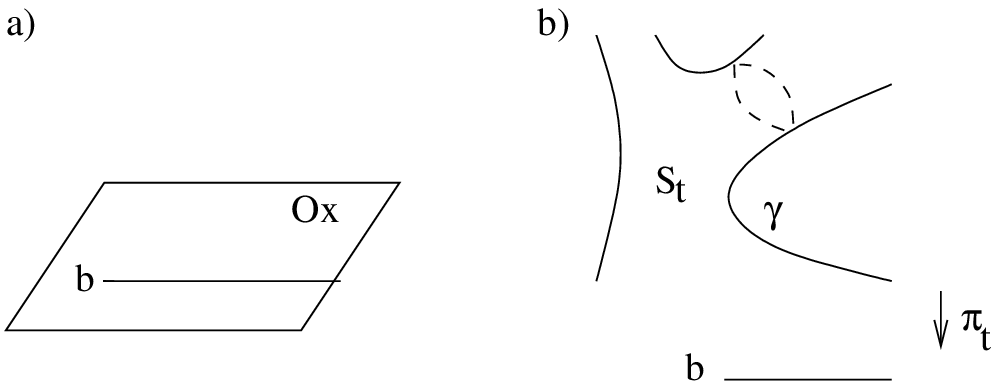}
    \caption{}
    \label{fig:2}
  \end{center}
\end{figure}

\begin{proof} A generic ray issued from $b$ has a pair of liftings
to $S_t$ under the projection, whose oriented union $\gamma$
connects two distinct points at infinity of the compactified level curve
(the orientations of the liftings are opposite). It follows from
the previous Proposition that some generator has a nonzero intersection index
with $\gamma$. Hence, the projection of any its representative intersects
the ray. For nongeneric ray the same statement follows by passing to the limit.
The Corollary is proved.
\end{proof}

\def\ddo{D_{\delta_0}}

We prove Theorem \ref{t3.5} by contradiction. Suppose the contrary: let there be
a \crp polynomial $H$ whose all critical values are contained in the disc
$D_{\delta_0}$. Fix a noncritical value $t\in D_{\delta_0}$. Lemma \ref{lmon} implies that the vanishing cycles along all the
paths $\a\subset\ddo$ from $t$ to 0
generate the integer homology group of $S_t$
(the latter paths, whose corresponding vanishing cycles generate
the homology, may be chosen critically-regular by Remark \ref{rkcreg}).
We show that in fact,
they cannot generate it (due to smallness of $\delta_0$). The contradiction thus
obtained will prove Theorem \ref{t3.5}. To do this, we prove that there exist a critical
value $b\in Ox$ of the projection $\pi_t$ and a disc $D_r\subset Ox$ whose closure
is disjoint from $b$, $|b|>r$, such that the canonical representatives
of the previous vanishing cycles are all projected inside $D_r$. Hence, their
projections do not intersect the ray issued from $b$ away from 0. Thus, by Corollary
\ref{corint}, these cycles cannot generate the homology group.

By definition, the projection of the canonical representative of vanishing cycle
is a piecewise linear curve whose vertices are critical values of $\pi_t$.
We have to show that some of the latters lies outside some disc
$D_r$, in particular, is distant from 0. The first step to do this is the next Lemma,
which gives an a priori lower bound of the maximal distance between critical
values of $\pi_0$. This is the main technical Lemma of the Section.
Here the critical values are understood in the following
generalized sense.

\begin{definition} Let $H$ be a ultra-Morse polynomial.
A {\it generalized critical value} of the projection $\pi_t:S_t\to Ox$
is either a critical value of $\pi_t$ (at a critical point
where the curve $S_t$ is regular),
or the projection image of a critical point of $H$ in $S_t$.
For any $t$ the set of the generalized
critical values of the projection $\pi_t$ will be denoted by $CS_t$,
as in 1.4.
\end{definition}

\begin{remark} The projection image of a critical point of $H$ is
a double generalized critical value of projection, provided
that transversality condition (\ref{locbr}) holds. This follows from
Proposition \ref{doublecr}.
\end{remark}

\begin{proposition} \label{preta}
Let $H$ be an ultra-Morse polynomial, and
zero lines of  its highest homogeneous part be transversal to the
$y$- axis. Then for any $t\in\cc$ the number of generalized critical
values of $\pi_t$ (with multiplicities) is equal to
$$\eta(n)=n(n+1).$$
\end{proposition}
The Proposition follows from Bezout's theorem, as the similar statement
of Remark \ref{reta}.

\begin{lemma} \label{distcr0}
Let $H$ be a \crp polynomial of degree $n+1\geq3$,
$(x,y)$  be orthogonal coordinates
in $\ccd$ satisfying (\ref{disty0}). There exists a generalized critical value
$b\in CS_0$ such that
\begin{equation}
 |b|>r(n),\ r(n)=(c'(H))^{7n^2}n^{-35n^2}.\label{rn}\end{equation}
\end{lemma}

Lemma \ref{distcr0} is proved in 2.4.

Let $D_r$ be a circle that separates the previous value $b$ from 0, i.e., $|b|>r$,
and such that
\begin{equation} dist (\partial D_r, CS_0)\geq \frac{r(n)}{2\eta(n)},\
\eta(n)=n(n+1).\label{distdr}
\end{equation}
(Its existence follows from Proposition \ref{preta} and (\ref{rn}).)
We show that the disc $D_r$ is a one we are looking for. To do this, we prove the
next Lemma.

\begin{lemma} \label{ldr} Let $H$, $(x,y)$, $r(n)$ be as in the previous Lemma,
$\partial D_r$ be a circle satisfying (\ref{distdr}). Then the points from
$CS_t$ do not cross the circle, while $t$ ranges in $\ddo$.
\end{lemma}

Lemma \ref{ldr} is proved in 2.5, where we also prove the following
more general

\begin{lemma} \label{deltae} Let $H$ be a \crp polynomial of degree $n+1\geq3$, $(x,y)$ be orthogonal
coordinates in $\ccd$ that satisfy (\ref{disty0}). Fix arbitrary $x\in Ox\setminus CS_0$
and put
$$\var=\min(dist(x,CS_0),1).\ \text{Then}$$
\begin{equation} x\notin CS_t\ \text{for any}\ t\in D_{\Delta(n,\var)},
\ \Delta(n,\var)=(c'(H))^{4n^3}n^{-17n^3}\var^{n(n+1)}.
\label{de} \end{equation}
\end{lemma}

The proofs of Lemmas \ref{distcr0} and \ref{deltae} use a priori bounds
from 2.3 for unit-scaled polynomials.

\begin{proof} {\bf of Theorem \ref{t3.5}.} Let the orthogonal
coordinates in $\ccd$ be chosen to satisfy (\ref{disty0}) and
(\ref{locbr}). Recall that the
canonical representative of any cycle in $S_t$
vanishing to 0 along a path in $\ddo$ is
projected onto a piecewise linear curve with vertices in $CS_t$.
All the vertices lie in $D_r$ (and hence, the projection itself also by convexity).
For the local vanishing
cycle this statement follows from definition. The same statement for the global vanishing cycle then follows from Lemma \ref{ldr}, convexity
and continuity.
This together with the discussion at the beginning of the Subsection proves Theorem \ref{t3.5} modulo Lemmas \ref{distcr0} and \ref{ldr}.
\end{proof}

\subsection{From \crp to weakly-normalized. Proof of Theorem \ref{topxy} and
Addendum modulo Lemmas \ref{distcr0} and \ref{ldr}.}

\def\la{\lambda}

Here we deduce Theorem \ref{topxy} from Theorem \ref{t3.5}.

Let $H$ be a weakly-normalized polynomial: then $||h||_{max}=1$, and the critical values
of $H$ lie in the disc $\overline D_2$. Consider the auxiliary polynomial
\begin{equation}
\wt H=\la^{-(n+1)} (H(\la x, \la y)-H(0)),\ \la>0, \
\text{such that}\
||\wt H'||_{max}=1.
\label{wth}\end{equation}

The possibility of choice of such $\la$ follows from the condition saying that $H'-H(0)\not\equiv0$, see the beginning of the paper.
By construction, the new polynomial $\wt H$ is \crp and
has the same highest homogeneous part $h$. For the proof of Theorem \ref{topxy} we prove
the following upper bound of $\la$ using Theorem \ref{t3.5}.

\begin{proposition} Let $\la$ be  as in (\ref{wth}), $\delta_0$ be
as in Theorem \ref{t3.5}. Then
\begin{equation}\la\leq\la_0=(\frac4{\delta_0})^{\frac1{n+1}}.\label{lao}
\end{equation}
\end{proposition}

\begin{proof} Let $a_i$, $\wt a_i=\la^{-(n+1)}(a_i-H(0))$ be the critical values of $H$ and $\wt H$ respectively. By weak
normalizedness, $|a_i|\leq2$, and $H(0)$ is one
of the $a_i$' s (in particular, $|H(0)|\leq2$).
This together with the previous formula implies that
$|\wt a_i|\leq4\la^{-(n+1)}$. On the other hand, there exists an $i$
such that $|\wt a_i|\geq\delta_0$ (Theorem \ref{t3.5}). The two latter inequalities imply
(\ref{lao}).
\end{proof}

Denote $\wt H_{\wt\la}$ a polynomial given by formula (\ref{wth})
(with $\la$ replaced by arbitrary $\wt\la$).
Then $||\wt H_{\la_0}'||_{max}\leq1$ (the previous Proposition
and the monotonicity of the norm $||\wt H_{\wt\la}'||_{max}$ as
a function in $\wt\la$). Thus, the polynomial $\wt H_{\la_0}$ is
unit-scaled.

Let $t\in\cc$, $|t|\leq5$, $S_t=\{ H=t\}$, $\wt S_{\tau}=\{\wt H_{\la_0}=\tau\}$, $\tau=\la_0^{-(n+1)}(t-H(0))$. By definition,
$$S_t=\la_0\wt S_{\tau}, \ |\tau|<5 \ \text{(formula (\ref{lao})
and inequality}\ |H(0)|\leq2).$$
By Proposition \ref{phl} applied to $\wt H_{\la_0}$, the topology of the
curve $\wt S_{\tau}$ lies in the bidisc $D_{X_n,Y_n}$, see (\ref{xnyn}). This together
with the previous formula implies that the topology of $S_t$ lies in the bidisc
$D_{X,Y}$, $X=\la_0X_n$, $Y=\la_0Y_n.$ Now elementary inequalities
imply that $X<Y<R_0$. This proves Theorem \ref{topxy}.

Now let us prove the Addendum to Theorem \ref{topxy}.
We have proved above that
$||\wt H_{\la_0}'||_{max}\leq1$, $\la_0=(\frac4{\delta_0})^{\frac1{n+1}}$. On the other hand,
it follows from definition that $||\wt H'_{\la_0}||_{max}\geq\la_0^{-(n+1)}
(||H'||_{max}-|H(0)|)$. This together with the previous inequality
and the fact that $|H(0)|\leq2$ and (\ref{fdel}) yield
$$||H'||_{max}\leq\la_0^{n+1}+2=4(c'(H))^{-13n^4}n^{63n^4}+2<
(c'(H))^{-13n^4}n^{64n^4}.$$
This proves the Addendum.

\subsection{A priori bounds of topology of unit-scaled polynomials}

In the present Subsection we prove the following more precise
version of Proposition \ref{phl}.

\begin{lemma} \label{topxn} Let $H(x,y)=h+H'$ be a unit-scaled
polynomial (see Definition \ref{dunitsc})
of degree $n+1\geq3$. Let the orthogonal coordinates in
$\mathbb C^2$ satisfy (\ref{disty0}).
Then for any $t\in\mathbb C$ with $|t|\leq5$ the topology of
the level curve $S_t=\{ H=t\}$ is contained in
the bidisc $D_{X_n,Y_n}$ (see Definition \ref{deftop}),
where $X_n$, $Y_n$ are the same, as in (\ref{xnyn}).
Moreover, the complement $S_t\setminus D_{X_n,Y_n}$ is a union of
graphs $y=y_i(x)$ of $n+1$ functions $y_0(x),\dots,y_n(x)$
holomorphic in $\mathbb C\setminus D_{X_n}$ such that
\begin{equation} |y_i(x)-y_j(x)|>\frac{c'(H)}{3n}|x| \
\text{for any}\ x\in\mathbb C\setminus D_{X_n},\ i\neq j.\label{yijy}\end{equation}
\end{lemma}

\begin{proof} Consider the
scalar product $(\ ,\ )$ on the two-dimensional vectors in $\mathbb C^2$
defined by the standard Hermitian metric. The
$y$- axis is not a zero line of $h$ by (\ref{disty0}), so, the latter
may be written as
\begin{equation}h(x,y)=c_{-1}\prod_{i=0}^n(y-c_ix)=c_{-1}
\prod_{i=0}^n(E, v_i),\ v_i=(-\bar c_i,1),\ E=(x,y)\ \text{is the
Euler field}.\label{prodh}\end{equation}

The zero lines of the homogeneous polynomial $h(x,y)$ are $y=c_ix$.
As it is shown below, inequality (\ref{disty0}) and
unit-scaledness condition imply the following a priori
estimates of the $c_i$' s:
\begin{equation} (n+1)^{-\frac{n+1}2} < c_{-1}\leq1, \label{c-1}
\end{equation}
\begin{equation} |c_i|<\sqrt n\ \text{for any}\ i=0,\dots,n.
\label{ci}\end{equation}
\begin{equation}|c_i-c_j|>\frac{c'(H)}n\ \text{for any}\ i,j\geq0,\ i\neq j.
\label{cij}
\end{equation}

Afterwards, we prove the statement of the Lemma as follows. Fix an
$x$, $|x|\geq X_n$, and consider the polynomial $h(x,y)$ as that with fixed
$x$ and variable $y$. The polynomial $h$ has $n+1$ roots
$c_ix$, $i\geq0$. The distances between them are bounded from below by
$\frac{c'(H)}n|x|$, see (\ref{cij}).
Using this bound we show that if we add
to $h$ a polynomial $H'$ of smaller degree and at most unit norm and then
subtract a $t$, $|t|\leq5$, then the new polynomial
$P=h+H'-t$ has a root $y_i(x)$ in the $\delta$- neighborhood of each $c_ix$,
where
\begin{equation}\delta=\frac{c'(H)}{3n}|x|.\label{xdelta}\end{equation}
(The proof of this statement is based on the next
more general Proposition \ref{root}.)
The previous neighborhoods are disjoint by (\ref{cij}), hence,
the roots $y_i(x)$ are distinct and thus, depend holomorphically on
$x\notin D_{X_n}$. Inequality (\ref{yijy}) of the Lemma
follows immediately from (\ref{cij}) and (\ref{xdelta}).

The Lemma says that the topology of the curve $S_t$, $|t|\leq5$, is contained
in $D_{X_n,Y_n}$. To prove this, now we have  to show that
$|y||_{\{(x,y)\in S_t,\ x\in D_{X_n}\}}<Y_n$. To do this, it suffices to prove
the inequality
\begin{equation}|y_i|(x)< Y_n,\ \text{whenever} \ |x|=X_n, \ i=0,\dots,n.
\label{yi<<yn}\end{equation}
We will prove it at the end of the Subsection.
Then the multivalued extensions of the $y_i$' s to the interior of $D_{X_n}$
satisfy the same inequality by the maximum principle.
This proves the Lemma.

\begin{proof} {\bf of (\ref{ci})}.
 The distance of a zero line
$l_i=\{ y=c_ix\}$ of $h$ to the $y$- axis is greater than $\frac1{\sqrt n}$
by (\ref{disty0}).
On the other hand, this distance is equal to $\arctan{\frac1{|c_i|}}$. This
follows from the same statement in the case, when $c_i$ is real (one
can make an individual $c_i$ real by applying a rotation in the
coordinate $x$; the rotation preserves distances between lines). Therefore,
\begin{equation}\arctan{\frac1{|c_i|}}>\frac1{\sqrt n},\  \text{hence}\
|c_i|<\frac1{|\tan{\frac1{\sqrt n}}|}.\label{tanci}\end{equation}
This together with the classical inequality
$\tan x>x$ implies (\ref{ci}).
\end{proof}

\begin{proof} {\bf of (\ref{c-1})}. By unit-scaledness,
$||h||_{max}=\max_{|x|^2+|y|^2=1}|h|(x,y)=1.$
The previous maximum of $|h|$ is no less than $|h|(0,1)=|c_{-1}|$. Hence,
$|c_{-1}|\leq1$, which proves the right inequality in (\ref{c-1}).

Now let us prove the lower bound of $c_{-1}$ from (\ref{c-1}).
Let $v_i$ be the vectors from the expression (\ref{prodh}) for $h$.
By definition and the same formula (\ref{prodh}),
\begin{equation}1=||h||_{max}\leq|c_{-1}|\prod_{i=0}^n||v_i||.
\label{hmac}\end{equation}
By definition and (\ref{ci}),
$||v_i||=\sqrt{1+|c_i|^2}<\sqrt{n+1}.$ Substituting the
latter inequality to (\ref{hmac}) yields
$1<|c_{-1}|(n+1)^{\frac{n+1}2}$. This proves (\ref{c-1}).
\end{proof}

\begin{proof} {\bf of (\ref{cij})}. Fix $i\neq j\geq0$.
Consider the line $x=1$ and its segment bounded
by its intersections with the lines $y=c_ix$ and $y=c_jx$.
By definition, the length of this segment is equal to
$|c_i-c_j|$ and
is greater than the angle between the two latter lines
(since the previous length is no less than
the tangens of the angle and by the inequality
$\tan x > x$). The latter angle
is equal to the distance between the lines, and hence, is
no less than
$\frac{c'(H)}n$ by definition. This proves (\ref{cij}).
\end{proof}

\begin{proposition} \label{root}
Let $p(y)=p_0(y-c_0)\dots(y-c_n)$ be a polynomial
of degree $n+1$ in one variable, $q(y)$ be a polynomial of a smaller degree,
$P(y)=p(y)+q(y)$. Let $\delta>0$ be such that
\begin{equation}|c_i-c_j|>2\delta\ \text{for any}\ i\neq j,
\label{cijd}\end{equation}
\begin{equation}\max_{|y|\leq\max_i|c_i|+\delta}|q(y)|<|p_0|\delta^{n+1}.
\label{maxyc}\end{equation}
Then the polynomial $P(y)=p(y)+q(y)$ has a root in the $\delta$- neighborhood
of each root $c_i$ of the polynomial $p(y)$.
\end{proposition}
\begin{proof} The statement of the Proposition holds true for the
initial polynomial $p(y)$ with the roots $c_i$. To prove it for $P(y)$,
we consider the auxiliary family of polynomials $P_s(y)=p(y)+sq(y)$,
$s\in[0,1]$. Let us show that for any $s\in[0,1]$ the polynomial $P_s(y)$
does not vanish on the circles $|y-c_i|=\delta$
(whose closed discs are disjoint by (\ref{cijd})). This will prove that
the roots of $P_s$ (which depend continuously on $s$) will  not leave
the $\delta$- neighborhoods of $c_i$, as $s$ is changed from 0 to 1.
It suffices to show that
\begin{equation}|p(y)|>|q(y)|, \ \text{whenever}\ min_i|y-c_i|=\delta.
\label{pqy}\end{equation}
The right-hand side $|q(y)|$ is less than $|p_0|\delta^{n+1}$ by
(\ref{maxyc}). The left-hand side is no less than
$|p_0|\prod|y-c_i|\geq|p_0|\delta^{n+1}$, since
$|y-c_i|\geq\delta$. This proves (\ref{pqy}) and the Proposition.
\end{proof}

\begin{proof} {\bf of (\ref{yijy}).}
Fix an $x$, $|x|\geq X_n$. Let $\delta$ be the same, as in (\ref{xdelta}). The distance between the zeros of the polynomial
$p(y)=h(x,y)$ is no less than
$3\delta$ by (\ref{cij}) and (\ref{xdelta}). Put $q(y)=H'(x,y)-t$.
To show that the polynomial $P(y)=H(x,y)-t=p(y)+q(y)$ has a root in
the $\delta$- neighborhood of each $c_ix$,
let us apply Proposition \ref{root} (here $p_0=c_{-1}$).
To check the conditions of the Proposition, it suffices to show that
\begin{equation}|H'(x,y)-t|<|c_{-1}|\delta^{n+1}\ \text{whenever}\
 |y|\leq|x|\max_{i\geq0}|c_i|+\delta,\ |t|\leq5.\label{previn}\end{equation}
 To prove the latter inequality, let us estimate its left-hand
side (which is a polynomial of degree no greater than $n$) in terms of
$||H'||_{max}$. It follows from definition that
\begin{equation}
|H'(x,y)|\leq||H'||_{max}(|x|^2+|y|^2)^{\frac n2}\leq
||H'||_{max}(|x|+|y|)^n, \ \text{whenever}\ |x|^2+|y|^2\geq1,
\label{hmaxy}\end{equation}
which is the case, since $|x|\geq X_n$.
Substituting inequality $||H'||_{max}\leq1$
(unit-scaledness), the second inequality in (\ref{previn}),
$\delta=\frac{c'(H)}{3n}|x|$ yields
$$|H'(x,y)|\leq (|x|(1+\max_i|c_i|+\frac{c'(H)}{3n}))^n<
(2|x|\sqrt n)^n$$
by (\ref{ci}) and the inequality $c'(H)\leq1$. Hence,
$|H'(x,y)-t|\leq2(2|x|\sqrt n)^n$, since $|t|\leq5$ and
$|x|\geq X_n$.

Now for the proof of
(\ref{previn}) it suffices to show that for $|x|\geq X_n$
$$2(2|x|\sqrt n)^n<|c_{-1}|\delta^{n+1}=
|c_{-1}|(\frac{c'(H)}{3n})^{n+1}|x|^{n+1}, \ \text{i.e.,}$$
\begin{equation} |x|>2^{n+1}(\sqrt n)^n(3n)^{n+1}|c_{-1}|^{-1}
(c'(H))^{-(n+1)}.
\label{xn>?}\end{equation}
The latter right-hand side is less than $X_n$ by
(\ref{c-1}) and elementary inequalities. This proves
(\ref{xn>?}) and (\ref{previn}). Inequality (\ref{yijy})
is proved.
\end{proof}

\begin{proof} {\bf of inequality (\ref{yi<<yn}).} For any
$x$ with $|x|=X_n$ one has
$|y_i(x)-c_ix|<\delta=\frac{c'}{3n}X_n$, as it was proved before. This
together with (\ref{ci}) and elementary inequalities yields
(\ref{yi<<yn}). The proof of Lemma \ref{topxn} is completed.
\end{proof}
\end{proof}

\subsection{Lower bound of the maximal critical value of projection.
Proof of Lemma \ref{distcr0}}
Let $H$ be a \crp polynomial. For simplicity denote
\begin{equation}r=r(n)=(c')^{7n^2}n^{-35n^2}, \ \text{see (\ref{rn}).}
\label{rrn}\end{equation}
Lemma \ref{distcr0} says that the projection $\pi:S_{0}\to Ox$ has a
generalized critical value with module greater than
$r$. To prove this, we have to show that the
functions $y_i(x)$, which give the roots of $H(x,y)$ (as a polynomial in $y$ with fixed $x$),
cannot have global holomorphic branches with disjoint graphs outside the disc $\overline{D_r}$. (A point of intersection of graphs
is a critical point of $H$. Hence, its $x$- coordinate is a
generalized critical value of the projection.)

Let us firstly sketch a proof of
this statement in the simplest case, when the polynomial
$H(x,y)$ under consideration is a product of $n+1$ linear (nonhomogeneous)
functions: $H(x,y)=c_{-1}\prod_{i=0}^n(y-c_ix-b_i)$.
We know that $H(0)=0$, so, at least one $b_i$ is zero, say, $b_0$. On the
other hand, $H$ has nonzero nonconstant lower terms,
moreover, their max- norm
is unit. This implies that at least one $b_i$ is not zero and the maximal
module of the $b_i$'s admits a lower bound (see (\ref{b1}) below,
let $b_1$ have the maximal module).
The line $y=c_1x+b_1$ intersects the line $y=c_0x+b_0=c_0x$ at a
point with the $x$- coordinate $x_0=\frac{b_1}{c_0-c_1}$.
An explicit calculation using (\ref{rrn}),
(\ref{cij}) and the lower bound for $|b_1|$
shows that  $|x_0|>r$. Therefore, the graphs of the linear functions
$y_i(x)=c_ix+b_i$ are not disjoint outside the disc $\overline{D_r}$.

Now let us prove the statement of the Lemma in the general case,
when the functions $y_i(x)$ are not necessarily linear. If they
are not holomorphic
outside $\overline D_r$, then the Lemma follows immediately.
Now suppose $y_i(x)$ are holomorphic outside $\overline D_r$. We
show (using smallness of $r$) that that two of them
(say, $y_0$ and $y_1$) have intersected graphs over
the complement to $\overline D_r$. This together with the discussion
at the beginning of the Subsection proves the Lemma.

By the previous holomorphicity assumption,
\begin{equation}y_i(x)=c_ix+b_i+\phi_i(x),\ \phi_i(x)\to0,\ \text{as}\
x\to\infty,\ \phi_i(x)\ \text{is holomorphic outside}\
\overline{D_r}.\label{yphi}\end{equation}
Let $X_n,$ $Y_n$ be the constants from (\ref{xnyn}).
Using Lemma \ref{topxn}, (\ref{yphi}) and smallness of the radius $r$,
we show that the functions $\phi_i(x)$ are small on a domain distant from $D_r$:
\begin{equation}|\phi_i(x)|< \frac{3Y_nr}{|x|-r} \ \text{outside} \
\overline{D_r},
\label{phi<}
\end{equation}
so, the functions $y_i(x)$ are close to the linear functions
$c_ix+b_i$ and the polynomial $H$ is close
to the product $c_{-1}\prod_{i\geq0}(y-c_ix-b_i)$.
As it will be shown below,
this together with previous arguments
in the case of product of linear functions imply that
one of the $b_i$'s, say $b_0$, is small, and the other one (say, $b_1$) is
large:
\begin{equation}\text{there exists a}\ b_i\ \text{(say,}\ b_0) \ \text{with}\
|b_0|<
6Y_nr^{\frac1{n+1}},\label{b0}\end{equation}
\begin{equation}\text{there exists a}\ b_i\ \text{(say,}\ b_1) \ \text{with}\
|b_1|>
\frac1{8n^2Y_n},\label{b1}
\end{equation}
We prove that the graphs of $y_0(x)$ and $y_1(x)$ are
intersected over a point $x\notin\overline{D_r}$.

 The graphs of the linear functions are intersected at a point with the $x$-
 coordinate
\begin{equation}x_0=\frac{b_1-b_0}{c_0-c_1}.\ \text{We show that}\ |x_0|>4r
\label{x04r}\end{equation}
(in fact, this follows from (\ref{cij}), (\ref{b0}), (\ref{b1})).
Then we show that the closeness of the $y_i(x)$' s to the linear functions
implies that the difference $y_0(x)-y_1(x)$ vanishes at some point
in the $\frac{|x_0|}2$- neighborhood of $x_0$. The latter neighborhood
(and hence, the latter point) are disjoint from $\overline{D_r}$ by
(\ref{x04r}). This proves the Lemma.

In the proof of the inequalities mentioned above we use the following
a priori bounds of $y_i$ and $b_i$. The multivalued functions $y_i$
satisfy the bound
\begin{equation}|y_i(x)|\leq Y_n,\ \text{whenever}\ |x|\leq X_n,\
i=0,\dots,n\label{yi<yn}
\end{equation}
(Lemma \ref{topxn}, which says that the topology of $S_0$ is contained in
$D_{X_n,Y_n}$). By the residue formula and (\ref{yphi}),
$$b_i=\frac1{2\pi i}\int_{|x|=X_n}\frac{y_i(x)}x. \ \text{Therefore,}$$
\begin{equation}|b_i|\leq Y_n.\label{bi<y}\end{equation}

\begin{proof} {\bf of (\ref{phi<}).} The maximal
value on $\overline{D_{X_n}}$ of the (multivalued) function $\phi_i$ is less
than $3Y_n$:
$$|\phi_i(x)|=|y_i(x)-c_ix-b_i|\leq Y_n+|c_i|X_n+|b_i|\leq2Y_n+\sqrt n X_n<
3Y_n$$
(by (\ref{yi<yn}), (\ref{bi<y}), (\ref{ci}) and (\ref{xnyn})).
The function $\phi_i$ is holomorphic in $|x|>r$ and vanishes at infinity.
Hence, the Cauchy formula implies that for any $x\notin\overline{D_r}$ one has
$$\phi_i(x)=\frac1{2\pi i}\int_{|\zeta|=r}
\frac{\phi_i(\zeta)}{x-\zeta},\ \text{so},\
|\phi_i(x)|\leq\frac{r\max_{|\zeta|=r}|\phi_i(\zeta)|}{|x|-r}.$$
This together with the previous inequality implies (\ref{phi<}).
\end{proof}

\begin{proof} {\bf of (\ref{b0}).} By assumption,
$H(x,y)=c_{-1}\prod_{i=0}^n(y-y_i(x))=c_{-1}\prod_{i=0}^n
(y-c_ix-b_i-\phi_i(x)),\
H(0,0)=0$. Therefore, up to a polynomial vanishing at 0, $H$ is equal to
\begin{equation}(-1)^{n+1}c_{-1}
(\prod_{i=0}^nb_i+\sum_i(\phi_i(x)\prod_{j\neq i}y_j(x))).\label{hbi}
\end{equation}
(In fact, the sum in (\ref{hbi}) is a polynomial in $x$.)
The polynomial (\ref{hbi}) vanishes at 0, as does $H$.
Therefore, the module
of the constant term $\prod b_i$ should be no greater than the
maximal module
of the sum in (\ref{hbi}) on any circle centered at 0, e.g., on
$\partial D_1\subset\mathbb C\setminus\overline{D_r}\cap D_{X_n}$.
The latter maximal module is no greater than
$$(n+1)\max_i
\max_{\partial D_1}|\phi_i|Y_n^n\leq\frac{3(n+1)Y_n^{n+1}r}{1-r}<
4r(n+1)Y_n^{n+1}\ \text{(see (\ref{phi<})). Hence,} $$
$$\prod_{i=0}^n|b_i|<4r(n+1)Y_n^{n+1},\ \text{so},\
\min_i|b_i|<(4r)^{\frac1{n+1}}(n+1)^{\frac1{n+1}}Y_n<6r^{\frac1{n+1}}Y_n
$$
\end{proof}

\begin{proof} {\bf of (\ref{b1}).} By definition,
$H(x,y)=c_{-1}\prod_{i=0}^n(y-c_ix-b_i-\phi_i(x))$, hence,
\begin{equation}
H'(x,y)=-\sum_{i=0}^n(b_i+\phi_i(x))\frac{H(x,y)}{y-y_i(x)}.
\label{h'bi}\end{equation}
As it is shown below, this together with the equality
$||H'||_{max}=1$ implies that the $b_i$' s cannot be too small
simultaneously. To do this,
we use the following relation between the max- norm of a
(nonhomogeneous) polynomial and the maximum of its module on the unit bidisc:
\begin{equation}\text{for any polynomial}\ G(x,y),\ deg G\leq n, \
\text{one has}\
||G||_{max}\leq 2n\max_{|x|, |y|\leq1}|G(x,y)|.\label{gxy}\end{equation}
\begin{proof}  Let
\begin{equation}G(x,y)=\sum g_{ij}x^iy^j.\ \text{Then}\
\max_{|x|, |y|\leq1}|G(x,y)|\geq\sqrt{\sum|g_{ij}|^2}\label{gxy2}\end{equation}
Indeed, the polynomial $G$ can be considered as a
Fourier polynomial in the harmonics of
the torus $\{|x|,|y|=1\}$, the
polynomial coefficients coincide with the Fourier coefficients. The maximum
of the module of the polynomial on the bidisc is equal to that
(of the
Fourier polynomial) on the torus. This follows from the
one-dimensional maximum principle applied to $G$ along the
lines $x=const$, $y=const$. The maximal module of the Fourier
 polynomial   is no less than its  $L_2$- norm (divided by the area
of the torus). This
proves (\ref{gxy2}). On the other hand, it follows from definition that
\begin{equation}||G||_{max}\leq\sum|g_{ij}|.
\label{gdef}\end{equation}
The total number
of monomials of degree no greater than $n$ is equal to $\frac{(n+1)(n+2)}2$.
Now the classical mean inequality says that
$$\sum|g_{ij}|\leq\sqrt{\frac{(n+1)(n+2)}2\sum|g_{ij}|^2}.$$
This together with (\ref{gxy2}), (\ref{gdef}) and elementary
inequalities implies
(\ref{gxy}).
\end{proof}

Now it follows from (\ref{gxy}) that
\begin{equation}\max_{|x|,|y|\leq1}|H'|\geq\frac1{2n}||H'||_{max}=
\frac1{2n}.\label{maxmax}\end{equation}
On the other hand,
$$\max_{|x|,|y|\leq1}|H'|=\max_{|x|=1,|y|\leq1}|H'|
\ \text{(maximum principle)}.$$
Now substituting to (\ref{h'bi}) the upper bounds
(\ref{phi<}), (\ref{yi<yn}) of $\phi_i|_{|x|=1}$ and $y_i$ respectively
(and also the inequality $c_{-1}\leq1$, see (\ref{c-1})),
yields
\begin{equation}
\max_{|x|,|y|\leq1}|H'|=\max_{|x|=1,|y|\leq1}|H'|\leq
(n+1)(\max_i|b_i|+\frac{3Y_nr}{1-r})\max_i\max_{|x|=1,|y|\leq1}
|\frac{H(x,y)}{y-y_i(x)}|.\label{ynry}\end{equation}
Let us firstly prove that for any $i$ and any $x$, $|x|=1$, one has
\begin{equation}\max_{|y|\leq1}
|\frac{H(x,y)}{y-y_i(x)}|\leq 10^{\frac{n+1}2}.
\label{fracy}\end{equation}
Fix an $x$, $|x|=1$, and consider the unit disc
$D_1(y_0(x))\subset Oy$ centered at $y_0(x)$.
If it is disjoint from the unit disc centered at 0, then the previous
fraction module maximum is no greater than the
module maximum of $|H|$ on the same unit disc (by definition).
Otherwise, the same statement
holds but for the maximums evaluated on the disc $D_3\supset
(D_1(y_0(x))\cup D_1(0))$, since then $dist(\partial D_3,y_0(x))>1$.
This together with the maximum principle
imply that in both cases the fraction maximum from (\ref{fracy})
is less than $\max_{|y|\leq3}|H(x,y)|\leq\max_{|x|\leq1,|y|\leq3}
|H(x,y)|$. The latter maximum is no greater than
$$||H||_{max}\max_{|x|\leq1,|y|\leq3}(|x|^2+|y|^2)^{\frac{n+1}2}=
10^{\frac{n+1}2},$$
as in (\ref{hmaxy}). This proves (\ref{fracy})

Now substituting (\ref{fracy}), (\ref{maxmax})
and the (elementary) inequality $\frac{3Y_nr}{1-r}<4Y_nr$
to (\ref{ynry}) yields
$$\frac1{2n}<2n10^{\frac{n+1}2}(\max_i|b_i|+4Y_nr),\ \text{thus,}\
\max_i|b_i|>\frac{1-16n^210^{\frac{n+1}2}Y_nr}{4n^210^{\frac{n+1}2}}.$$
This together with elementary inequalities implies (\ref{b1}).
\end{proof}

\begin{proof} {\bf of (\ref{x04r}).}
Thus, by (\ref{b0}), (\ref{b1}) we have
$|b_0|<6Y_nr^{\frac1{n+1}}<\frac1{8n^2Y_n}<|b_1|.$
The intermediate inequality holds and moreover,
its third term is at least twice greater
than its second term (by elementary inequalities). Hence,
$|b_1|>2|b_0|, \ \text{thus,}\ |b_1-b_0|>\frac{|b_1|}2.$
Recall that
$|c_1-c_0|\leq|c_1|+|c_0|<2\sqrt n$ by (\ref{ci}). This together with
(\ref{b1}) and the previous inequality implies (\ref{x04r}):
\begin{equation}|x_0|=\frac{|b_1-b_0|}{|c_0-c_1|}>\frac{|b_1|}{2|c_0-c_1|}>
\frac1{32n^2Y_n\sqrt n}=(32n^{\frac52}Y_n)^{-1}>4r.\label{x0b}
\end{equation}
\end{proof}

Now for the proof of Lemma \ref{distcr0} it suffices to
show that there exists a point $x$ in the disc $D_{\frac{|x_0|}2}(x_0)$
such that $y_0(x)=y_1(x)$. To prove this,
we use the argument principle: we consider the difference
$y_1(x)-y_0(x)$ restricted to the circle $|x-x_0|=\frac{|x_0|}2$ and show that
it vanishes nowhere and the increment of its argument along the circle
is $2\pi$. This is true for $y_i$ replaced by the linear functions
$c_ix+b_i=y_i(x)-\phi_i(x)$: the difference of the latters
is equal to $(c_1-c_0)(x-x_0)$. In the general case one has
$$y_1(x)-y_0(x)=(c_1-c_0)(x-x_0)+(\phi_1(x)-\phi_0(x)).$$
The restriction to the circle of the difference of the linear functions
(the first term in the right-hand side of the previous formula)
has a constant module $|c_1-c_0|\frac{|x_0|}2$. To show that
the difference of the $y_i$' s has the same argument increment, as the previous
first term, it suffices to prove that the second term
(the difference of the $\phi_i$'s) has a smaller module:
\begin{equation}|\phi_1(x)-\phi_0(x)|<|c_1-c_0|\frac{|x_0|}2\
\text{whenever}\ |x-x_0|=\frac{|x_0|}2.\label{x0c}\end{equation}
For the proof of the latter inequality let us estimate
its right-hand side.
By (\ref{cij}), $|c_1-c_0|>\frac{c'(H)}n$, hence,
\begin{equation}|c_1-c_0|\frac{|x_0|}2>\frac{c'(H)}{2n}|x_0|.
\label{c01}\end{equation}
Now let us estimate from above the difference of the $\phi_i$'s
on the circle  $|x-x_0|=\frac{|x_0|}2$: by (\ref{phi<}) and inequality
$r<\frac{|x_0|}4$ (see (\ref{x04r})),
$$|\phi_1(x)-\phi_0(x)|\leq|\phi_0(x)|+|\phi_1(x)|<
2\frac{3Y_nr}{|x|-r}\leq\frac{6Y_nr}{\frac{|x_0|}2-r}<\frac{24Y_nr}{|x_0|}.$$
By (\ref{c01}), for the proof of (\ref{x0c}) it suffices to check that
$\frac{24Y_nr}{|x_0|}<\frac{c'}n\frac{|x_0|}2$, or equivalently,
$$r<c'(48nY_n)^{-1}|x_0|^2.$$
Indeed, by (\ref{x0b}), the right-hand side of the latter inequality is greater than
$$c'(48nY_n)^{-1}(32n^{\frac52}Y_n)^{-2}>
c'(64nY_n)^{-1}(32n^{\frac52}Y_n)^{-2}=
c'2^{-16}n^{-6}Y_n^{-3}\geq c'n^{-22}Y_n^{-3}.$$
Recall that $Y_n=n^{8n}(c')^{-2n}$ by definition,
hence, the previous right-hand side is equal to
$(c')^{6n+1}n^{-24n-22}>(c')^{7n}n^{-35n}>r$.
Inequality (\ref{x0c}) is proved. The proof of Lemma \ref{distcr0} is completed.

\subsection{Noncritical values of projection of variable
level curve. Proof of Lemmas \ref{ldr} and \ref{deltae}}
Firstly we prove Lemma \ref{deltae} and then deduce Lemma
\ref{ldr}.

Let $X_n$, $Y_n$ be the constants from (\ref{xnyn}).
Let $y_i(x,t)$, $i=0,\dots,n$, be the roots of the polynomial $H(x,y)$ in $y$,
which are multivalued holomorphic functions
with branching at those points $(x,t)$ where $x\in CS_t$.
The conditions of Lemma \ref{topxn} hold ($H$ is unit-scaled), hence,
the topology of $S_t$, $|t|\leq5$, is contained in $D_{X_n,Y_n}$. This means
that for any $t$, $|t|\leq5$,

- the functions $y_i(x,t)$ are holomorphic
in $x\notin\overline{D_{X_n}}$ and have disjoint graphs;

- for any branch of $y_i$ one has $|y_i|(x)\leq Y_n$, whenever
$|x|\leq X_n$.

The first statement implies that  the set $CS_t$ of generalized
critical values lies in $\overline D_{X_n}$,
whenever $|t|\leq5$.

Now let us take any $x\in D_{X_n}\setminus CS_0$ (recall that
$\var=\min(dist(x,CS_0),1)$).
Let us prove that then $x\notin CS_t$ whenever $|t|<\Delta(n,\var)$
(by definition, $\Delta(n,\var)<5$).
To do this, consider the discriminant (taken up to constant)
of the polynomial $H$ in $y$:
$$\Sigma_t(x)=\prod_{i<j}(y_i(x,t)-y_j(x,t))^2.$$
The zeros of $\Sigma_t$ coincide
with the generalized critical values of the projection of $S_t$.
For any fixed $t\in\mathbb C$ $y_i(x,t)=c_ix(1+o(1))$, as $x\to\infty$.
Therefore (denote $x_l(t)$ the roots of $\Sigma_t$, $l=1,\dots,n(n+1)$),
\begin{equation} \Sigma_t(x)=\prod_{0\leq i<j\leq n}(c_j-c_i)^2\prod_{l=1}^{n(n+1)}
(x-x_l(t)).\label{sigmacprod}\end{equation}
Let $|t|<\Delta(n,\var)$ (in particular, $|t|<4$). Let us show that
$x$ is not a zero of $\Sigma_t$. To do this, we
use the following a priori lower bound of $\Sigma_0(x)$:
\begin{equation}|\Sigma_0(x)|\geq(\prod_{0\leq i<j\leq n}|c_j-c_i|^2)
\var^{n(n+1)}>
(\frac{c'(H)}n\var)^{n(n+1)}=(c'(H))^{n(n+1)}n^{-n(n+1)}\var^{n(n+1)}.\label{sigma0>}\end{equation}
This  follows from inequalities $|x-x_i|\geq\var$ (which hold by assumption),
(\ref{sigmacprod}) and (\ref{cij}). To show that $\Sigma_t(x)\neq0$ for small
$t$, we estimate from above the derivative of $\Sigma_t(x)$ in $t$ by using
the following
a priori upper bound of $\Sigma_t(x)$ valid whenever $|t|\leq5$ and $|x|\leq X_n$:
$$|\Sigma_t(x)|=\prod_{i<j}|y_i-y_j|^2(x,t)\leq
(2Y_n)^{n(n+1)}\leq(nY_n)^{n(n+1)}$$
\begin{equation}=n^{8n^2(n+1)+n(n+1)}(c')^{-2n^2(n+1)}
<n^{16n^3}(c')^{-3n^3}\label{sig<}\end{equation}
(since
$|y_i(x,t)|\leq Y_n$ and by elementary inequalities).

The classical inequality on derivative of holomorphic function
(which follows from Cauchy lemma) says that if $\psi(t)$ is a bounded
function holomorphic in $t\in D_R$, $|\psi|<\mathcal M$, then
$|\psi'(t)|\leq\frac{\mathcal M}{R-|t|}$. Applying this inequality
to the function
$\Sigma_t(x)$, now with fixed $x$ and variable $t\in D_4$, $R=5$,
$\mathcal M=n^{16n^3}(c')^{-3n^3}$ (see (\ref{sig<})),
yields that for any $t$ with $|t|<\Delta(n,\var)<4$
$$|(\Sigma_t(x))'_t|\leq\frac{\mathcal M}{5-|t|}\leq\mathcal M=n^{16n^3}
(c')^{-3n^3}, \ \text{hence, by (\ref{de}),}$$
$$|\Sigma_t(x)-\Sigma_0(x)|\leq\Delta(n,\var)n^{16n^3}(c')^{-3n^3}
=(c')^{4n^3}n^{-17n^3}\var^{n(n+1)}n^{16n^3}(c')^{-3n^3}$$
\begin{equation}=
(c')^{n^3}n^{-n^3}\var^{n(n+1)}\label{sigmat0}\end{equation}
Now to prove that $\Sigma_t(x)\neq0$ it suffices to show that the
 right-hand side of (\ref{sigmat0}) is less than $\Sigma_0(x)$. But it is
 clearly less than the right-hand side in (\ref{sigma0>}), which gives a lower
 bound of $|\Sigma_0(x)|$.
This proves the inequality $\Sigma_t(x)\neq0$ and Lemma \ref{deltae}.

\begin{proof} {\bf of Lemma \ref{ldr}.} To show that $CS_t$ does not
meet $\partial D_r$ whenever $|t|\leq\delta_0$,
we apply Lemma \ref{deltae} to each point $x\in\partial D_r$.
Then by (\ref{distdr})
$$\var=dist(x,CS_0)\geq\var'=\frac{r(n)}{2\eta(n)},\ \eta(n)=n(n+1).$$
We show that
\begin{equation}\Delta(n,\var')>\delta_0\label{ddoo}
\end{equation}
 This together with the inequality
$\Delta(n,\var)\geq\Delta(n,\var')$ (monotonicity)
implies that $\Delta(n,\var)>\delta_0$. Together with
Lemma \ref{deltae}, this yields $x\notin CS_t$ whenever
$|x|=r$, $|t|\leq \delta_0$ and proves Lemma \ref{ldr}.

By (\ref{de}),
$$\Delta(n,\var')=(c')^{4n^3}n^{-17n^3}(\frac{r(n)}{2n(n+1)})^{n(n+1)}=
(c')^{4n^3}n^{-17n^3}(r(n))^{n(n+1)}(2n(n+1))^{-n(n+1)}.$$
This together with
formula (\ref{rn}) for $r(n)$ and elementary inequalities imply
that the latter right-hand side (and hence, the $\Delta$)
is greater than
$$(c')^{4n^3}n^{-17n^3}(c')^{11n^4}n^{-53n^4}n^{-3n^3}\geq(c')^{13n^4}n^{-63n^4}=
\delta_0.$$
Lemma \ref{ldr} is proved.
\end{proof}

\section{Lower bounds of the formula for the main determinant.
Proof of Theorem \ref{tlowerch}}

Here we prove Theorem \ref{tlowerch} (in 3.6), which gives a lower bound of
$C(h,\Omega)$. To do this, we prove lower bounds of the terms
of its formula (\ref{2.21}). In 3.1-3.4
we prove Lemma \ref{lpd>} (lower bound of $P_d$). In 3.5 we prove
an upper bound of the discriminant $\Sigma(h)$:
\begin{equation} \Sigma(h)<n^{6n^2}.\label{sigma>}\end{equation}
 In 3.7 we prove a lower bound of the
constant $C_n$:
\begin{equation}C_n>e^{-12n^2}.\label{cn>}\end{equation}

\subsection{Lower bound of $P_d$. The sketch of the proof of Lemma \ref{lpd>}}

By definition, $P_d=\det A_d(h,\Omega(d))$,
where $A_d$ is the $(d+1)\times(d+1)$-
matrix defined in 1.8. Its lines are naturally identified
with vectors in the space of complex polynomials of degree $d+1$,
which are split into two collections:

- the collection (denoted by $\Pi$) of $2(d-n+1)$ vectors
$x^jy^{d-n-j}\frac{\d h}{\d y}$, $x^jy^{d-n-j}\frac{\d h}{\d x}$, which
do not depend on the forms $\o_i$;

- the collection (denoted $\Pi_{\Omega(d)}$) of $2n-d-1$ vectors
$\frac{D\o_i}{n-i+1}$,
$\o_i=x^{l'(i)}y^{m'(i)+1}dx$, $l'(i)+m'(i)=d$, which
depend only on $\Omega(d)$: by definition,
\begin{equation}D\o_i=(m'(i)+1)x^{l'(i)}y^{m'(i)}. \label{doilm}
\end{equation}

We have to prove a lower bound of the maximal value of
$P_d$ as a function of variable monomial form tuple $\Omega(d)$.
Firstly we prove its next a priori lower bound in terms of
complex volume of the collection $\Pi$, see the following Definition.

\begin{definition} Let $\Pi=\{ v_1,\dots,v_k\}\in\mathbb C^m$
be  a collection of vectors, $k\leq m$. Consider
the real parallelogramm formed by the vectors $v_j$
and $iv_j$. The  {\it complex volume} of $\Pi$
 (denoted $Vol\Pi$) is the square root of the real $2k$- dimensional
standard Hermitian (Euclidean) volume of the previous parallelogramm.
\end{definition}

\begin{remark} \label{vol=det}
The complex volume is nonzero, if and only if
the vector collection is linearly independent over complex numbers.
If $k=m$, then the complex volume is equal to the module of the
determinant of the $m\times m$- matrix whose lines are formed by
the components of the collection vectors.
\end{remark}

\begin{proposition} \label{pchoice} Let $h$ be a generic
homogeneous polynomial of degree $n+1\geq3$, $d\in\{n,\dots,2n-2\}$,
$\Pi$ be the corresponding
vector collection from the beginning of the Subsection.
There exists a tuple $\Omega(d)$ of $2n-d-1$ forms of the type
$x^ly^{m+1}dx$, $l+m=d$, such that
\begin{equation} P_d(h,\Omega(d))>n^{-4n}Vol\Pi.\label{pi'}
\end{equation}
\end{proposition}
The Proposition is proved at the end of the Subsection by elementary
linear algebra arguments.

The principal part of the Section is the proof of lower bound
of $Vol\Pi$: we show that if $||h||_{max}=1$, then
\begin{equation}Vol\Pi>n^{-42n^2}(c'(h))^{6n^2}\label{vol}\end{equation}
This together with Proposition \ref{pchoice} implies
Lemma \ref{lpd>}.

To prove (\ref{vol}), we consider the following space of two-dimensional
vector polynomials:
\begin{equation}V_s=\{ v=(v_1(x,y),v_2(x,y)),\ deg v_i=s,\
v_i\ \text{are homogeneous}\}, \ s\leq2n-2, \end{equation}
equipped with the standard Hermitian structure.
Denote $||v||_2$ the corresponding Hermitian norm. The space $V_s$ has
the standard orthonormal
basis of $2(s+1)$ monomials $(x^iy^j,0)$, $(0,x^iy^j)$, $i+j=s$.

Consider the linear operator
\begin{equation}L:V_s\to V_{n+s},\ L(v)=\frac{d h}{dv}=v_1
\frac{\partial h}{\partial x}+v_2\frac{\partial h}{\partial y}
\label{lgrad}\end{equation}
Let $s\leq n-2$, $d=n+s$. By definition,
the vectors of the collection $\Pi$ are the images under $L$
of the previous basic monomials. Denote
$$\nu(L)=||L^{-1}|_{LV_s}||^{-1}=\min\{\frac{||Lv||_2}{||v||_2},\
v\in V_s\setminus0\}.$$
It follows from definition that $Vol\Pi\geq (\nu(L))^{2(d-n+1)}$.

In Subsection 3.4 we show that $\nu(L)>n^{-21n}(c'(h))^{3n}$,
or equivalently,
\begin{equation}||Lv||_2>n^{-21n}(c'(h))^{3n}\ \text{for each}\ v\in V_s\ \text{with}\
||v||_2=1,\ s\leq n-2.\label{lv>}\end{equation}
\begin{proof} {\bf of (\ref{vol}).}
The two previous inequalities imply that
$$Vol\Pi>(n^{-21n}(c'(h))^{3n})^{2(d-n+1)}\geq(n^{-21n}(c'(h))^{3n})^{2n}=
n^{-42n^2}(c'(h))^{6n^2}.$$
This proves (\ref{vol}) modulo (\ref{lv>}).
\end{proof}

To prove (\ref{lv>}), we consider the following extension of
the  max-norm to the vector polynomials:
\begin{equation}||v||_{max}=\max_{|x|^2+|y|^2=1}
\sqrt{|v_1|^2+|v_2|^2}(x,y),
\ v=(v_1,v_2)\in V_s.\label{extmax}\end{equation}

\begin{remark} \label{rinvar}
The max- norm of a polynomial or of a polynomial
vector field in $\ccd$ is invariant under orthogonal transformations,
while the Hermitian norm isn't.
\end{remark}

In 3.3 we prove an inequality similar to (\ref{lv>}), but for
the $max$- norm instead of the Hermitian norm:
\begin{equation}||Lv||_{max}>(c'(h))^{3n}n^{-19n}||v||_{max}.
\label{max>}\end{equation}
This is the main technical inequality of the Section.
Firstly we prove it for constant vector fields $v$.
To prove it for a higher degree homogeneous polynomial vector field $v$,
we choose appropriate orthogonal coordinates $(x,y)$
(see Remark \ref{rinvar}) and
consider the space of vector polynomials in $V_s$ proportional to
the Euler vector field with polynomial coefficient:
$$V_{s,e}=\{ v=(xQ(x,y),yQ(x,y))\in V_s\ |\  Q\ \text{is a homogeneous
polynomial of degree}\ s-1\}.$$
We decompose $v$ as a sum
\begin{equation}v=v'+v'',\ v'=(xQ,yQ)\in V_{s,e}, \ v''\perp V_{s,e}.
\ \text{One has}\ Lv=\frac{dh}{dv'}+\frac{dh}{dv''}.\label{v'v''}
\end{equation}
We estimate from below the contributions of $v'$ and $v''$ to
$Lv$ separately. Firstly we show that
\begin{equation}||\frac{dh}{dv'}||_{max}\geq n^{-4n}
||v'||_{max}.\label{dhv'>}
\end{equation}
Put
$$S=\{|x|^2+|y|^2=1\}.$$
On the zero lines of $h$ one has $\frac{dh}{dv'}=0$ (by definition), hence,
$\frac{dh}{dv}=\frac{dh}{dv''}$. We show that
\begin{equation}\max|\frac{dh}{dv''}||_{S\cap \{ h=0\}}\geq
(c'(h))^{3n}n^{-12n}||v''||_{max}
\ \text{for any} \ v''\in V_s,\ v''\perp V_{s,e}.
\label{dhv''>}\end{equation}
This together with the previous statement implies that
\begin{equation}||Lv||_{max}\geq(c')^{3n}n^{-12n}
||v''||_{max}.\label{dhv>1}
\end{equation}
We show that either the latter right-hand side already gives
itself the desired lower bound (\ref{max>})
of $||Lv||_{max}$, or the contribution of $v''$ to $Lv$
is dominated (many times) by that of $v'$. In the latter case
lower bound (\ref{max>}) will be deduced from
(\ref{dhv'>}) and the following
simple a priori bound (proved at the end of 3.3):
\begin{equation}||\frac{dh}{dw}||_{max}\leq2^{n+1}||w||_{max}
 \ \text{for any}\ w\in V_s, s<n, \ \text{whenever}\
||h||_{max}=1. \label{dhu<}\end{equation}

In the proofs of the previously mentioned lower bounds of
the derivatives of $h$ along the vector fields $v$, $v'$, $v''$
we use simple a priori bounds of coefficients of a polynomial in
terms of its $max$- norm and
the following relation between the Hermitian and the max- norms
(these bounds and relation will be proved in 3.2):
\begin{equation}\text{for any}\ s\in\mathbb N, Q\in V_s, \
\text{one has}\ \frac1{\sqrt{2(s+1)}}||Q||_{max}\leq||Q||_2\leq
2^{\frac s2}||Q||_{max}
\label{hermax},\end{equation}
in the case, when $Q$ is a scalar polynomial of degree $s$,
the left coefficient
$\frac1{\sqrt{2(s+1)}}$ may be replaced by $\frac1{\sqrt{s+1}}$.

\begin{proof} {\bf of Proposition \ref{pchoice}.}
Consider the collection $\Pi_{\Omega(d)}$ of the
$2n-d-1$ lines (defined by $\Omega(d)$) of the matrix $A_d$.
By definition, all the elements of the i-th line of the collection
are zeros except for the one standing in the column numerated by
the monomial $x^{l'(i)}y^{m'(i)}$. It follows from definition
and (\ref{doilm}) that the latter element equals
$\frac{m'(i)+1}{n-i+1}\geq\frac1n$. These elements form
a unique nonzero minor (denoted $M$)
of maximal size in the collection $\Pi_{\Omega(d)}$,
$|M|\geq n^{d+1-2n}>n^{-2n}$
by the previous inequality. The determinant
$P_d=\det A_d$ is thus equal
to $M$ times the complementary minor (denoted $M'$), which is a
maximal size minor in the collection $\Pi$. This together with
the previous inequality implies that
\begin{equation}P_d\geq n^{-2n}M'.\label{pdm'}\end{equation}
It follows from definition and Remark \ref{vol=det},
that $M'$ is equal to the complex volume
of the orthogonal projection of $\Pi$ along
the coordinate plane generated by the monomials $x^{l'(i)}y^{m'(i)}$.
We can choose the tuple $\Omega(d)$ so that
the latter monomials be arbitrarily given $2n-d-1$ distinct
monomials of degree $d$ with unit coefficient,
or equivalently, the complementary coordinate plane orthogonal to them
be arbitrarily given coordinate plane of dimension $2(d-n+1)$.
Let us choose $\Omega(d)$ so that the previous orthogonal projection
of $\Pi$ have the maximal possible complex volume $M'$.

The triangle inequality says that $Vol\Pi$ is no greater than the
sum of the complex volumes of the orthogonal projections of $\Pi$ to
all the complex coordinate $2(d-n+1)$- planes in the space of
degree $d$ polynomials (which has dimension $d+1$). The number
of all the latter planes is equal to
$$C_{d+1}^{2(d-n+1)}=
\frac{(2n-d)\dots(d+1)}{(2(d-n+1))!}<(2n-d)^{2(d-n+1)}<n^{2n}.$$
Therefore, the maximal complex volume $M'$ of projection is no less than
\newline
$(C^{2(d-n+1)}_{d+1})^{-1}Vol\Pi>n^{-2n}Vol\Pi$. This
together with (\ref{pdm'}) implies (\ref{pi'}).
\end{proof}

\subsection{The max- norm and the coefficients. A priori bounds}

In the proof of (\ref{lv>}), (\ref{max>}) and (\ref{hermax}) we use the following properties of the max- norm.

\begin{proposition} \label{ppmaxho} Let $g(x,y)$ be a complex homogeneous
polynomial of degree $k$,
$$g=g_0\prod_{i=1}^k((x,y),u_i),\ u_i\in\mathbb C^2,\ ||u_i||=1.
\ \text{Then}$$
\begin{equation}||g||_{max}\leq|g_0|\leq||g||_{max}(2\sqrt k)^k.\label{maxho}
\end{equation}
(Here the expression $((x,y),u_i)$ is the standard Hermitian scalar product
(linear in $(x,y)$) of $u_i$ and the Euler vector field $(x,y)$.
\end{proposition}
\begin{proof} In the proof of Proposition \ref{ppmaxho} we use the following relation
between distance of complex lines and the Hermitian scalar product.

\begin{proposition} \label{pdistscal}
Let $u_1,u\in\mathbb C^2$, $||u_1||=||u||=1$. Let $l_1=u_1^{\perp}$ be
the complex line orthogonal to $u_1$, $\lambda$ be the complex line containing $u$. Then
the Hermitian scalar product $(u_1,u)$ admits the lower bound
\begin{equation}|(u_1,u)|\geq\frac12dist(l_1,\lambda).\label{distscal}
\end{equation}
\end{proposition}
\begin{proof} It follows from definition that
$|(u_1,u)|$ is equal to the length of the orthogonal
projection of $u$ along the line $l_1$. The latter length is equal to the sine
of the minimal angle between $u$ and a real line in $l_1$. The latter
angle is equal to $dist(l_1,\lambda)\leq\frac{\pi}2$ by definition, thus,
$$|(u_1,u)|=\sin{dist(l_1,\lambda)}\geq\frac12dist(l_1,\lambda). \ \text{This proves
(\ref{distscal}).}$$
\end{proof}

The left inequality in (\ref{maxho}) follows from definition.
Let us prove the right inequality. Denote $l_i=u_i^{\perp}$ the zero
line of $g$ orthogonal to $u_i$. By Proposition \ref{pdisty0},
 there is a complex
line $\lambda$ through the origin such that $dist(\lambda, l_i)>\frac1{\sqrt k}$ for
all $i$ (let us fix such a line $\lambda$ and a unit vector
$u$ on it). Then by (\ref{distscal}),
$$|g(u)|=|g_0|\prod|(u,u_i)|,\ |(u,u_i)|\geq\frac12dist(\lambda,l_i)>
\frac1{2\sqrt k}. \ \text{Hence},$$
$$||g||_{max}\geq|g(u)|>|g_0|(\frac1{2\sqrt k})^k.$$
\end{proof}

{\bf The max- norm of product:}
\begin{equation}||gQ||_{max}\geq (\frac1{2\sqrt{k+m}})^{k+m} ||g||_{max}||Q||_{max}\
\text{for any}\ k,m\in\mathbb N\label{maxpr}\end{equation}
$$\text{and homogeneous polynomials} \
g,Q, \ k=degg, m=degQ.$$

\begin{proof} Let us write
down the decompositions preceding (\ref{maxho})
of the polynomials $g$, $Q$, $gQ$
(denote $g_0$, $Q_0$, $(gQ)_0$ the corresponding constant
factors from their decompositions). By definition and the left inequality
in (\ref{maxho}), $|(gQ)_0|=|g_0||Q_0|\geq||g||_{max}||Q||_{max}$.
Applying the right inequality in (\ref{maxho}) to $gQ$ yields
$$||gQ||_{max}>|(gQ)_0|(\frac1{2\sqrt{k+m}})^{k+m}.$$
 This together with the previous inequality proves (\ref{maxpr}).
\end{proof}

{\bf Hermitian and the max- norms: proof of (\ref{hermax}).} We use the inequality
\begin{equation}
\max_{|x|,|y|\leq1}|Q(x,y)|\leq\max_{|x|,|y|\leq 1}
(|x|^2+|y|^2)^{\frac s2}||Q||_{max}=2^{\frac s2}||Q||_{max},\
deg Q=s,\label{discball}\end{equation}
which follows from definition, as (\ref{hmaxy}).
The right inequality in (\ref{hermax}) follows from (\ref{discball}) and
(\ref{gxy2}) (inequality (\ref{gxy2}) holds true both for scalar and vector
polynomials).
Let us prove the left one. It follows from definition that $||Q||_{max}$ is no greater
than the sum of moduli of the coefficients of $Q$. The polynomial $Q$ is homogeneous
of degree $s$, hence, the number of its coefficients is $s+1$ (in the case of scalar
polynomial) and $2(s+1)$ in the vector case.
By the mean inequality, the latter sum is no greater than the sum of the squared moduli
of the coefficients times the square root of their number (either of $s+1$ or of
$2(s+1)$). This proves the left inequality in (\ref{hermax}).

\subsection{Lower bound of $||Lv||_{max}$. Proof of (\ref{max>})}
We choose the orthogonal coordinates $(x,y)$ as in (\ref{disty0}):
the distance of the $y$- axis to each zero
line of $h$ is greater than $\frac1{\sqrt n}$.

Let us firstly prove (\ref{max>}) in the case,
when $v$ is a constant vector polynomial, i.e.,
$s=deg v_1=deg v_2=0$ (we suppose that $||v||_2=1$, denote $l$ the complex line
containing $v$). Let  $S=\{|x|^2+|y|^2=1\}$.
By definition, the $\frac{c'(h)}{2n}$- neighborhoods in $S$
of the $n+1\geq3$ zero lines of the polynomial $h$ are disjoint. Therefore, the previous neighborhood of at
least one zero line (denote the latter zero line by $l_i$)
is disjoint from the line $l$.
Let us choose a point $z\in S\cap l_i$ and show (below) that
\begin{equation}|\frac{dh}{dv}|(z)>(c'(h))^{n+1}n^{-4n}.
\label{dh>}\end{equation}
The latter module of derivative is no greater than
$||\frac{dh}{dv}||_{max}=||Lv||_{max}$). This implies
(\ref{max>}) in the particular case under consideration modulo
(\ref{dh>}). Let us prove the  latter.

Let $u_i$ be the vector from the expression for $h$ in (\ref{maxho})
that is orthogonal to the
chosen zero line $l_i$: $(z,u_i)=0$. Then by Proposition \ref{ppmaxho},
\begin{equation}\frac{dh}{dv}(z)=h_0(v,u_i)\prod_{j\neq i}(z,u_j),
|h_0|\geq||h||_{max}=1.\label{dhomax}
\end{equation}
Let us estimate the factors in the latter right-hand side.
By Proposition \ref{pdistscal},
$$|(v,u_i)|\geq\frac12dist(l_i,l), \ |(z,u_j)|\geq\frac12dist(l_j,l_i).$$
The
former distance is no less than $\frac{c'(h)}{2n}$, and the latter one
is greater than $\frac{c'(h)}n$ by definition. Hence, by (\ref{dhomax}),
$$|\frac{dh}{dv}(z)|\geq\frac{c'}{4n}(\frac{c'}{2n})^n=
(c')^{n+1}2^{-n-2}n^{-n-1}\geq(c')^{n+1}n^{-n-2-n-1}>
(c')^{n+1}n^{-4n}.$$
This proves (\ref{dh>}) and hence, (\ref{max>}) in the case of constant vector field $v$.

Now let us prove (\ref{max>}) in the case, when $1\leq s=deg v_i
\leq n-2$.

\begin{proof} {\bf of (\ref{dhv'>}).} By homogeneity,
\begin{equation}
Lv'=\frac{dh}{dv'}=Q(x,y)(x\frac{\d h}{\d x}+y\frac{\d h}{\d y})=(n+1)hQ.\label{eul}\end{equation}

Applying inequality (\ref{maxpr})
yields that $||Lv'||_{max}=(n+1)||hQ||_{max}$ is no less than
$(\frac1{2\sqrt{n+s}})^{n+s}||h||_{max}||Q||_{max}$. By definition,
$||h||_{max}=1$, $||v'||_{max}=||Q||_{max}$, hence,
$$||Lv'||_{max}\geq(\frac1{2\sqrt{n+s}})^{n+s}||v'||_{max}
>(2\sqrt{2n})^{-2n}||v'||_{max}\geq n^{-4n}||v'||_{max}.$$
\end{proof}

\begin{proof} {\bf of (\ref{max>}) modulo (\ref{dhv''>}) and (\ref{dhu<}).}
By (\ref{dhv'>}), (\ref{dhu<}) and the triangle inequality,
\begin{equation}||\frac{dh}{dv}||_{max}\geq
||\frac{dh}{dv'}||_{max}-||\frac{dh}{dv''}||_{max}
\geq n^{-4n}||v'||_{max}-2^{n+1}||v''||_{max}
\label{dhvv-}\end{equation}

Let us consider the case, when
$||v''||_{max}<n^{-6n}||v'||_{max}$. Then the latter inequality implies that
the first term containing $||v'||_{max}$ in the right-hand side of
(\ref{dhvv-}) is at least twice greater than the second one.
By the same inequality,
$||v'||_{max}>\frac12||v||_{max}$. Therefore,
\begin{equation}||\frac{dh}{dv}||_{max}\geq\frac12n^{-4n}||v'||_{max}>
\frac14n^{-4n}||v||_{max}\geq
n^{-5n}||v||_{max}.\label{dhv>2}\end{equation}
This proves (\ref{max>}) in the case under consideration.

Now let us consider the opposite case, when
$||v''||_{max}\geq n^{-6n}||v'||_{max}$. Then
$$||v''||_{max}(1+n^{6n})\geq||v''||_{max}+||v'||_{max}
\geq||v||_{max}.$$
This together with (\ref{dhv>1}) (which follows from (\ref{dhv''>}))
  implies (\ref{max>}):
$$||\frac{dh}{dv}||_{max}\geq(c')^{3n}n^{-12n}
||v''||_{max}\geq(c')^{3n}n^{-12n}(1+n^{6n})^{-1}||v||_{max}>
(c')^{3n}n^{-19n}||v||_{max}.$$
\end{proof}

In the proof of lower bound  (\ref{dhv''>}) of the contribution
of $v''$ to $Lv$ we use the
 following description of the subspace of $V_s$ orthogonal to $V_{s,e}$:

a vector polynomial $v''\in V_s$ is
orthogonal to $V_{s,e}$, if and only if it has the form
\begin{equation}v''=r_1y^s\frac{\d}{\d x}+r_2x^s\frac{\d}{\d y}+R(x,y)
(x\frac{\d}{\d x}-y\frac{\d}{\d y}),\label{ortdec}\end{equation}
$$ r_1,r_2\in\cc,\ R\ \text{is a homogeneous polynomial of degree}\ s-1.$$
The statement saying that each vector polynomial from (\ref{ortdec})
is orthogonal to $V_{s,e}$
follows from definition. The inverse statement follows from the coincidence
of the dimensions of the space $V_{s,e}^{\perp}$ and the space of the vector polynomials
from (\ref{ortdec}): the former dimension equals $dim V_s-dim V_{s,e}=
2(s+1)-s=s+2$; the latter one equals $s+2$ by definition.

\begin{proof} {\bf of (\ref{dhv''>}).} We choose appropriate
zero line $l_i$ of $h$ (as follows) and show  that
\begin{equation}|\frac{dh}{dv''}||_{S\cap l_i}\geq
(c')^{3n}n^{-12n}||v''||_{max}.\label{dhvl>}\end{equation}
This will prove (\ref{dhv''>}).

Consider the lines through the origin that are
tangent to the vector field $v''$. They are defined by the equation
$yv_1''(x,y)-xv_2''(x,y)=0$ of degree $s+1$, so, their number is at
most $s+1\leq n$, which is less than the number
$n+1$ of the zero lines of $h$.
The $\frac{c'(h)}{2n}$- neighborhoods of the latters are
disjoint by definition. Hence, at least one zero line of $h$
(fix it and denote $l_i$) has distance no less than
$\frac{c'}{2n}$ from the lines tangent to $v''$.

Let us prove (\ref{dhvl>}). On the line $l_i$ one has
$x\frac{\d h}{\d x}+y\frac{\d h}{\d y}=(n+1)h=0$, hence,
$$\frac{d h}{d v''}|_{l_i}=v_2''\frac{\d h}{\d y}-
v_1''\frac yx\frac{\d h}{\d y}=(xv_2''-yv_1'')
x^{-1}\frac{\d h}{\d y}. \ \text{Therefore},$$
\begin{equation}|\frac{d h}{d v''}||_{S\cap l_i}
\geq |yv_1''-xv_2''||\frac{\d h}{\d y}||_{S\cap l_i}.
\label{v''xy}\end{equation}
Now we estimate from below the factors in the
right-hand side of (\ref{v''xy}). The second factor admits
the following lower bound:
 \begin{equation}|\frac{\partial h}{\partial y}||_{S
\cap l_i}>(c')^{n+1}n^{-4n}.\label{secfac}\end{equation}
This follows from (\ref{dh>}) applied to  $v=\frac{\d}{\d y}$.
(The condition of inequality (\ref{dh>})
saying that $dist(l,l_i)\geq\frac{c'}{2n}$ (in our case $l=Oy$)
is satisfied by choice of the coordinates, see the
beginning of the Subsection.)
Now let us estimate from below the first factor in the right-hand side
of (\ref{v''xy}): we show that
\begin{equation}|yv_1''-xv_2''|(z)\geq (c')^{2s}n^{-8s}||v''||_{max}\
\text{for any}\ z\in S\cap l_i.\label{area>}\end{equation}
This together with (\ref{v''xy}) and (\ref{secfac}) implies (\ref{dhvl>}):
$$|\frac{d h}{d v''}|(z)>(c')^{2s}n^{-8s}||v''||_{max}(c')^{n+1}n^{-4n}>
(c')^{3n}n^{-12n}||v''||_{max},\ \ z\in S\cap l_i.$$
For the proof of (\ref{area>}) let us firstly show that
\begin{equation}||yv_1''-xv_2''||_{max}\geq2^{-2s}||v''||_{max}
\ \text{for any}\ v''\perp V_{s,e}.\label{armax}\end{equation}
Then using the assumption that the line $l_i$ is distant from
the zero lines of the polynomial $yv_1''-xv_2''$ (which are tangent to
$v''$), we prove (\ref{area>}).

\begin{proof} {\bf of (\ref{armax}).} Firstly let us show that
$||yv_1''-xv_2''||_2\geq||v''||_2$.
A vector polynomial $v''\in V_{s,e}^{\perp}$
has the type (\ref{ortdec}), hence,
$$yv_1''-xv_2''=r_1y^{s+1}-r_2x^{s+1}+2xyR(x,y).$$
The three terms in the latter right-hand side are orthogonal to each other, as
are those in expression (\ref{ortdec}) for $v''$. Let us
compare the Hermitian norms of the corresponding terms in the latter
right-hand side and in (\ref{ortdec}).
The first (second) terms are equal.
The norm of the third one in the previous right-hand side is greater
than the norm of that in (\ref{ortdec}):
$$||2xyR(x,y)||_2=2||R||_2>\sqrt2||R||_2=
||R(x,y)(x\frac{\partial}{\partial x}-y\frac{\partial}{\partial y})
||_2.$$
This implies that $||yv_1''-xv_2''||_2\geq||v''||_2$.
Now applying (\ref{hermax}) twice yields
$$||yv_1''-xv_2''||_{max}\geq2^{-\frac{s+1}2}||yv_1''-xv_2''||_2\geq
2^{-\frac{s+1}2}||v''||_2\geq2^{-\frac{s+1}2}(2(s+1))^{-\frac12}||v''||_{max}.$$
The total coefficient at $||v''||_{max}$ in the latter right-hand side is no less than
$2^{-2s}$ (recall that we assume that $s\geq1$). This proves (\ref{armax}).
\end{proof}

\begin{proof} {\bf of (\ref{area>}).} Let $z\in S\cap l_i$. By (\ref{maxho}) and (\ref{armax}), one has
$$(yv_1''-xv_2'')(z)=q_0\prod_{j=0}^s(z,u_j),\ |q_0|\geq ||yv_1''-xv_2''||_{max}\geq
2^{-2s}||v''||_{max}.$$
Let $l_j'=u_j^{\perp}$ be a zero line of the polynomial $yv_1''-xv_2''$.
Recall that by the choice of $l_i$ one has
$dist(l_i,l_j')\geq\frac{c'}{2n}$
(see the beginning of proof of (\ref{dhv''>})). Hence, by Proposition \ref{pdistscal},
$|(z,u_j)|\geq\frac12dist(l_j',l_i)\geq\frac{c'}{4n}$.
This together with the previous formula and inequality implies that
$$|yv_1''-xv_2''|(z)\geq2^{-2s}||v''||_{max}(\frac{c'}{4n})^{s+1}=
2^{-4s-2}(c')^{s+1}n^{-s-1}||v''||_{max}$$
$$\geq(c')^{2s}n^{-8s}||v''||_{max}.$$
This proves (\ref{area>}). The proof of inequality
(\ref{dhv''>}) is completed.
\end{proof}
\end{proof}

\begin{proof} {\bf of (\ref{dhu<}) (upper bound of derivative).}
For any homogeneous vector polynomial $w$ of degree $s<n$
the module of the derivative $\frac{dh}{dw}$  is no greater than
\begin{equation}||\nabla h||\sqrt{|w_1|^2+|w_2|^2}\leq
||\nabla h||_{max}||w||_{max}=
 ||(\frac{\d h}{\d x},\frac{\d h}{\d y})||_{max}||w||_{max}.\label{gradhw}
 \end{equation}
Let us estimate from above the max- norm of the gradient of $h$, or
equivalently, the maximal
Hermitian norm of the gradient at points of the unit ball. To do this,
we use the well-known fact that the module of the
derivative of a holomorphic function (in one variable) at a given point is
no greater than the maximal module of the function in the unit disc centered
at the point under consideration. Therefore, the module of the derivative
of $h$ along any unit tangent vector (and hence, the Hermitian norm of the gradient)
at a point of the closed unit ball is no
greater than the maximal value of $h$ in the closed ball of radius 2.
By definition, $||h||_{max}=1$, hence,
$$|h(z)|\leq |z|^{n+1}\ \text{for any}\ z\in\mathbb C^2.$$
In particular,
\begin{equation}|h(z)|\leq 2^{n+1}\ \text{whenever}\ |z|\leq2.\label{h<2n}\end{equation}
Therefore, the max-norm of the gradient of $h$ is no greater than
$2^{n+1}$. This together with (\ref{gradhw}) implies (\ref{dhu<}).
\end{proof}

\subsection{From $max$- to Hermitian norm. Proof of (\ref{lv>})}
By (\ref{hermax}) and (\ref{max>}) (recall also that $s\leq n-2$),
$$||Lv||_2\geq\frac1{\sqrt{n+s+1}}||Lv||_{max}>
(2n)^{-\frac12}(c')^{3n}n^{-19n}||v||_{max}$$
$$\geq
n^{-2-19n}(c')^{3n}||v||_{max}\geq n^{-20n}(c')^{3n}
||v||_{max}.$$
By (\ref{hermax}), $||v||_{max}\geq 2^{-\frac s2}||v||_2>
2^{-\frac n2}||v||_2$, hence, the
right-hand side of the previous inequality is greater than
$$n^{-20n}2^{-\frac n2}(c')^{3n}||v||_2>n^{-21n}(c')^{3n}||v||_2.$$
This proves (\ref{lv>}).

\subsection{Upper bound of the discriminant. Proof of
(\ref{sigma>})} Consider a decomposition
\begin{equation}h(x,y)=\prod_{i=0}^n((x,y),v_i)\end{equation}
of the polynomial $h$ as a product of linear factors (not
necessarily of unit norm). The vectors $v_i$ are well-defined up to
multiplications by constants (the product of the latter constants should be unit).
Let $v_i=(v_{i,1},v_{i,2})$ be their components.
We use the following formula for the discriminant $\Sigma$ of $h$:
\begin{equation}\Sigma=\prod_{i< j}(v_{i,1}v_{j,2}-v_{i,2}v_{j,1})^2.
\label{sigmapr}\end{equation}
\begin{remark} The right-hand side of (\ref{sigmapr}) depends only
on $h$ and does not depend on the normalization of the vectors $v_i$.
\end{remark}
Consider a given pair of indices $i\neq j$ and
the corresponding pair of vectors $v_i$, $v_j$.
The module $|v_{i,1}v_{j,2}-v_{i,2}v_{j,1}|$
is equal to the complex volume of the vector collection $(v_i,v_j)$
(Remark \ref{vol=det}).  Therefore,
the previous module is no greater than $|v_i||v_j|$. Hence,
\begin{equation}|\Sigma|\leq\prod_{i<j}(|v_i||v_j|)^2=
\prod_{i\neq j}(|v_i||v_j|)=(\prod_i|v_i|)^{2n}.\label{sigmav<}
\end{equation}
Let us estimate the right-hand side of (\ref{sigmav<}).
To do this, consider the other  product decomposition
(preceding (\ref{maxho})) of the polynomial $h$.
Let $h_0=g_0$ be the corresponding coefficient from (\ref{maxho}).
By definition and (\ref{maxho}),
$$|h_0|=\prod_i|v_i|,\
|h_0|\leq ||h||_{max}(2\sqrt{n+1})^{n+1}=(2\sqrt{n+1})^{n+1}.$$
The two previous inequalities imply that
$$|\Sigma|\leq((2\sqrt{n+1})^{n+1})^{2n}=(2\sqrt{n+1})^{2n^2+2n}
<n^{2(2n^2+2n)}\leq n^{6n^2}.$$
This proves (\ref{sigma>}).

\subsection{Lower bound of $C(h,\Omega)$.
Proof of (\ref{lowerch})}

  Recall the formula for $C(h,\Omega)$ from 1.8:
$$C(h,\Omega)=C_n(\Sigma(h))^{\frac12-n}\prod_{d=n}^{2n-2}P_d,\
P_d=det A_d(h,\Omega).$$
Substituting inequalities (\ref{pd>}), (\ref{sigma>}) and
(\ref{cn>}) to its right-hand side yields
$$|C(h,\Omega)|> e^{-12n^2}n^{6n^2(\frac12-n)}(n^{-44n^2}(c')^{6n^2})^{n-1}
\geq e^{-12n^2}n^{-6n^3-44n^3}(c')^{6n^3}$$
$$=
n^{-12n^2(\ln n)^{-1}-50n^3}(c')^{6n^3}\geq
n^{-(6(\ln 2)^{-1}+50)n^3}(c')^{6n^3}>n^{-60n^3}(c')^{6n^3}.$$
Theorem \ref{tlowerch} is proved modulo (\ref{cn>}).

\subsection{Lower estimate of   the constant $C_n$. Proof of
(\ref{cn>})}

In the proof of (\ref{cn>}), we use the following inequalities:
\begin{equation}\ln2\pi>1;\label{6.19}\end{equation}
\begin{equation}N\ln N-N+1<\sum_{k=1}^N\ln k<(N+1)\ln(N+1)-N\   \text{for any}\ N\geq2;
\label{6.20}\end{equation}
\begin{equation}\sum_{k=a}^bk\ln k<(\frac{x^2\ln x}2-\frac{x^2}4)|_a^{b+1}
\ \text{for any}\ a<b,\ a,b\in\mathbb N.\label{6.21}\end{equation}
Inequality (\ref{6.20}) follows from the inequality
$$\int_1^N\ln xd x<\sum_{k=1}^N\ln k<\int_1^{N+1}\ln xd x$$
and the statement that the previous integrals are equal respectively to
the left- and right- hand sides of (\ref{6.20}). The previous   inequality follows
from increasing of the function $\ln x$. Inequality (\ref{6.21}) follows
from the statement that its left-hand side is less than
$\int_a^{b+1}x\ln xd x$ (increasing of the function $x\ln x$); the last
integral is equal to the right-hand side of (\ref{6.21}). By  (\ref{2.9}),
one has
\begin{equation}\ln|C_n|=\frac{n(n+1)}2\ln(2\pi)+\frac{n^2+n-4}2\ln(n+1)+
n\sum_{k=1}^{n+1}\ln k-\sum_{m=1}^{n-1}\sum_{k=1}^{m+n+1}\ln k.
\label{6.22}\end{equation}
Let us estimate the first three terms in the right-hand side of (\ref{6.22}).
The first term is greater than $\frac{n^2}2$ by (\ref{6.19}).   The second one is
greater than $\frac{n^2-2}2\ln n$. The third one (containing the sum till
$n+1$) is greater that $n((n+1)\ln(n+1)-n)>n^2\ln n-n^2$ by
(\ref{6.20}).
Substituting these
inequalities to (\ref{6.22}) and applying the   right inequality in (\ref{6.20}) to the
inner sum of the double sum in (\ref{6.22}) yields
$$\ln|C_n|>\frac{n^2}2+\frac{n^2-2}2\ln n+n^2\ln n-n^2-
\sum_{m=1}^{n-1}((m+n+2)\ln(m+n+2)-(m+n+1))$$
\begin{equation}=\frac{3n^2}2\ln n-\frac{n^2}2-
\ln n+\sum_{m=1}^{n-1}(m+n+1)-\sum_{k=n+3}^{2n+1}k\ln k.\label{6.23}
\end{equation}
Let us estimate the sums in the right-hand side of
(\ref{6.23}). The first sum is
equal to
\begin{equation}(n+1)(n-1)+\sum_{m=1}^{n-1}m=n^2-1+\frac{n(n-1)}2=
\frac{3n^2-n-2}2.\label{6.24}\end{equation}
The second sum is estimated by inequality (\ref{6.21}):
it is   less than
$$(\frac{x^2\ln x}2-\frac{x^2}4)|_{n+3}^{2(n+1)}$$
$$=\frac{4(n+1)^2(\ln(n+1)+
\ln2)}2-\frac{(n+3)^2\ln(n+3)}2-(n+1)^2+\frac{(n+3)^2}4$$
$$<\frac{3(n+1)^2\ln(n+1)}2+
2(n+1)^2\ln2-(n+1)^2+\frac{(n+3)^2}4$$
$$<\frac{3(n+1)^2\ln(n+1)}2+\frac{(n+3)^2}4+(n+1)^2.$$
Substituting the right-hand sides of (\ref{6.24}) and   the last inequality
instead of the first (respectively, second) sum in (\ref{6.23}) yields
$$\ln|C_n|>\frac{3n^2\ln n}2-\frac{n^2}2-\ln n+\frac{3n^2-n-2}2$$
\begin{equation}-
\frac{3(n+1)^2\ln(n+1)}2-(n+1)^2-\frac{(n+3)^2}4.\label{6.25}\end{equation}
Let us simplify inequality (\ref{6.25}). To do this, we use the following
inequalities:
\begin{equation}\ln(n+1)<\ln n+\frac1n;\ n+1<2n; \ n+3\leq\frac52n.
\label{6.26}\end{equation}
Let us estimate the fifth term in the right-hand side of
 (\ref{6.25}). By (\ref{6.26}),
it is less than
$$\frac{3(n+1)^2(\ln n+\frac1n)}2<\frac{3(n+1)^2\ln n}2+6n.$$
The sixth term in the same place is $(n+1)^2<4n^2$, and the seventh one is
less than $2n^2$ by (\ref{6.26}). Substituting these estimates to (\ref{6.25}) instead
of the corresponding terms yields
$$\ln|C_n|>\frac{3n^2\ln n}2-\frac{n^2}2-\ln n+\frac{3n^2-n-2}2-
\frac{3(n+1)^2\ln n}2-6n-4n^2-2n^2$$
$$=\frac32(n^2-(n+1)^2)\ln n+n^2-\frac{n+2}2
-\ln n-6n-6n^2$$
$$=-5n^2-\frac32(2n+1)\ln n-\ln n-\frac{13n+2}2.$$
Substituting the inequality $\ln n<\frac n2$ (valid for $n\geq2$)
to the right-hand side of the previous inequality yields
$$\ln|C_n|>-5n^2-\frac34(2n^2+n)-\frac n2-\frac{13n+2}2>-7n^2-9n>-12n^2.$$
This proves (\ref{cn>}).

\section{Acknowledgements}

I am grateful to Yu.S.Ilyashenko for helpful discussions and
inspiring major improvements of the text. I am also grateful
to V.Kharlamov for helpful discussions.
This paper was partly written during my  visits to the
Independent University of Moscow (this visit was supported by CNRS)
and State University of New York at Stony Brook. I wish to thank both Universities for hospitality and support. Research was supported
by part by RFBR grant 02-02-00482.

\end{document}